\documentclass{article} 
\usepackage{texdraw}\usepackage{graphicx}
\usepackage{amsmath} 
\usepackage{color}
\usepackage{colordvi}
\usepackage{texdraw}\usepackage{graphicx}
\usepackage{dcpic,pictex}

\def\newpic#1{%
   \def\emline##1##2##3##4##5##6{%
      \put(##1,##2){\special{em:point #1##3}}%
      \put(##4,##5){\special{em:point #1##6}}%
      \special{em:line #1##3,#1##6}}}
\newpic{}

\usepackage[cp1251,cp1257]{inputenc}\usepackage[english,latvian]{babel}\usepackage{amssymb,latexsym,amsfonts}
\newcommand{\apz}{\protect\makebox[1.8\width]{\rule[3.3pt]{8.5pt}{0.4pt} \hspace*{-12.6pt}\protect\raisebox{-.3\height}{$\leftharpoondown$}}}

\newtheorem{teor}{Theorem}[section]
\newtheorem{theorem}[teor]{Theorem}
\newtheorem{corollary}[teor]{Corollary}
\newtheorem{lemma}[teor]{Lemma}
\newtheorem{proposition}[teor]{Proposition}
\newtheorem{definition}[teor]{Definition}
\newtheorem{examples}[teor]{Examples}
\newtheorem{example}[teor]{Example}

\input cyracc.def
\font\tencyr=wncyr10
\def\cyr{\tencyr\cyracc}
\font\tencyi=wncyi10
\def\cyi{\tencyi\cyracc}

\begin{document}
\selectlanguage{english}
\title{Automaton (Semi)groups (Basic Concepts)}
\date{}
\author{J\=anis Buls, L\=iga U\v{z}ule, Aigars Valainis \\
{\small Department  of  Mathematics, University of Latvia, Ze\c{l}\c{l}u iela 25,}\\
{\small R\=\i ga, LV-1002 Latvia,
buls@fmf.lu.lv; liga.kulesa@lu.lv; avalains@gmail.com
}}
\def\keywords{\begin{center}{\bf Keywords}\end{center}
{automaton semigroups}}

\maketitle
\title{}
\abstract{\small{In this papper, we give an introduction to basic concepts of automaton semigroups. While we must note that this paper does not contain new results, it is focused on extended introduction in the subject and detailed examples.}}

\keywords{}

\section{ Preliminaries} 
\normalsize Let $A$ be a finite non-empty set and $A^*$ 
the free monoid generated by $A$. The set $A$ is also called an {\em 
alphabet}, its  elements are called {\em letters} and 
those of $A^*$ are called {\em finite words}. The identity element of $A^*$ is called 
an {\em empty word} and denoted by $\lambda$. We set 
$A^+=A^*\backslash\{\lambda\}$.

A word $w\in A^+$ can be written uniquely as a sequence of letters as 
$w=w_1w_2\ldots w_l$, with $w_i\in A$, $1\le i\le l$, $l>0$. The integer 
$l$ is called the {\em length} of $w$ and denoted by $|w|$. The length of 
$\lambda$ is 0. We set $w^0=\lambda$ and $\forall i \in \mathbb{N} \; w^{i+1}=w^iw\,.$ 

The word $w'\in A^*$ is a {\em factor} (or {\em subword}) of 
$w\in A^*$ if there exists $u,v\in A^*$ such that $w=uw'v$.
The words $u$ and $v$ are called, respectively, a {\em prefix} and a {\em suffix}. A pair $(u,v)$ is called an {\em occurrence} of $w'$ in $w$.
A factor $w'$ is called {\em proper} if $w\ne w'$. We denote, respectively, by F$(w)$, Pref$(w)$ and Suff$(w)$ the sets of $w$ factors, prefixes and suffixes.

An (indexed) infinite word $x$ on the alphabet $A$ is any total mapping 
$x\,:\,\mathbb{N}\rightarrow A$. We shall set for any $i\ge0$, $x_i=x(i)$
and write
\[
x=(x_i)=x_0x_1\ldots x_n\ldots \;.
\]
The set of all the infinite words over $A$ is denoted by $A^\omega$.

The word $w'\in A^*$ is a {\em factor} of $x\in A^\omega$ if there exists 
$u\in A^*$, $y\in A^\omega$ such that $x=uw'y$. 
The words $u$ and $y$ are called, respectively, a {\em prefix} and a {\em suffix}.
We denote, respectively, by F$(x)$, Pref$(x)$ and Suff$(x)$ the sets of $x$ factors, prefixes and suffixes. For any $0\le m\le n$, $x[m,n]$ denotes a factor $x_mx_{m+1}\ldots x_n$. The word 
$x[m,n]$ is called an {\em occurrence} of $w'$ in $x$ if $w'=x[m,n]$. The 
suffix $x_nx_{n+1}\ldots x_{n+i}\ldots$ is denoted by $x[n,\infty)$.

If $v\in A^+$, then we denote by $v^\omega$ the infinite word
\[
v^\omega=vv\ldots v\ldots \;.
\]

The {\em concatenation} of $u=u_1u_2\ldots u_k\in A^*$ and $x\in 
A^\omega$ is the infinite word 
\[
ux=u_1u_2\ldots u_kx_0x_1\ldots x_n\ldots
\]
For denoting concatenation we sometimes use symbol \#. 
The word $x$ is called {\em ultimately periodic} if there exists words $u\in A^*$, 
$v\in A^+$ such that
$x=uv^\omega$. In this case, $|u|$ and $|v|$ are called, respectively, an 
{\em anti-period} and a {\em period}. 

We use notation $\overline{0,n}$ to denote set $\left \{ 0,1,...,n \right \}$. 
\section{Serial composition of Mealy machines}\label{n5.2}

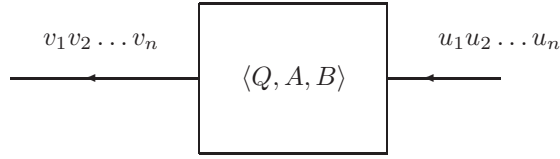
\begin{figure}[t]
\unitlength 1mm
\begin{picture}(100,30)(10,82)
\put(60,105){\line(0,-1){20}}
\put(60,85){\line(1,0){25}}
\put(85,85){\line(0,1){20}}
\put(85,105){\line(-1,0){25}}
\put(60,95){\line(-1,0){25}}
\put(100,95){\line(-1,0){15}}
\put(45.00,95.00){\vector(-1,0){0.2}}
\put(90.00,95.00){\vector(-1,0){0.2}}
\put(72.67,95.00){\makebox(0,0)[cc]{$\langle Q,A,B \rangle$}}
\put(100.00,100.00){\makebox(0,0)[cc]{$u_1u_2\ldots u_n$}}
\put(47.00,100.00){\makebox(0,0)[cc]{$v_1v_2\ldots v_n$}}
\end{picture}
\caption{An abstract Mealy machine. }
\label{z51}
\end{figure}

\begin{definition}
A 3-sorted algebra $V=\langle Q,A,B,\circ ,\ast \rangle$ is called  a 
{\em Mealy machine} if $Q,A,B$ are finite, nonempty sets, the mapping 
$Q\times A \stackrel{\circ}{\longrightarrow}Q$
is a total function and the mapping 
$Q\times A \stackrel{\ast }{\longrightarrow}B$
is a total surjective function. 
\end{definition}

The, set $Q$ is called {\em state set}, sets $A,B$ are called {\em input} and {\em output alphabet}, respectively.
The mappings $\circ$ and $\ast$ may be extended to $Q\times A^*$  by defining
\[
\begin{array}{lr}
q\circ \lambda =q,\quad & q\circ (ua)=(q\circ u)\circ a, \\
q\ast \lambda  =\lambda ,\quad & q\ast (ua)=(q\ast u)\#((q\circ u)\ast a)\,,
\end{array}
\] 
for each $q\in Q$, $(u,a)\in A^*\times A$. See \ref{z51}\,fig. for interpretation of Mealy machine as a word transducer. 
Henceforth, we shall omit parentheses if there is no danger of confusion. So, for example, we will write $q\circ u\ast a$ instead of $(q\circ u)\ast a.$ Similarly, we will write $q\circ q^{'}\ast a$ instead of $q\circ (q^{'}\ast a)$ where $q^{'} \in Q$.

\begin{figure}[b]
\hspace*{0.5cm} 
\centering
\scalebox{.9}{\input{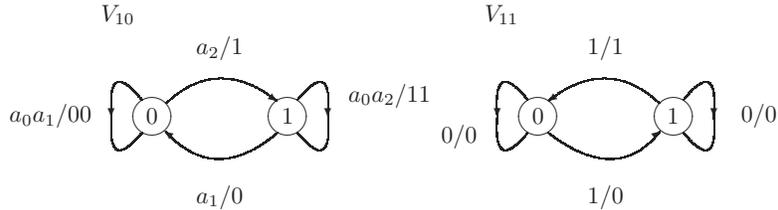}}
\caption{ Machines $V_{10}$ and $V_{11}$.}
\label{z52}
\end{figure}

\begin{definition}
A 3-sorted algebra $V_0=\langle Q,A,B,q_0,\circ ,\ast \rangle$ 
is called  an 
{\em initial Mealy machine} if $\langle Q,A,B,\circ ,\ast \rangle$ is a 
Mealy machine and $q_0\in Q$.
\end{definition}

Let $(q,x,y)\in Q\times A^\omega\times B^\omega$. 
We write $y=q*x$ if $\forall n \in \mathbb{N} \; y[0,n]=q*x[0,n]$ and say machine $V$ {\em transforms} $x$ to $y$. We say initial machine $V_0$ {\em transforms} $x$ to $y$ if $q=q_0$. We refer to words $x$ and $y$ as machines {\em input} and {\em output}, respectively. 

\begin{examples} 
\textnormal{
Look at \ref{z52}\ fig. for examples of machines $V_{10}$ and $V_{11}$. There we have \\
$V_{10}:a_0=(0,0),\; a_1=(0,1),\; a_2=(1,0)$; \\
$V_{11}:q\circ a=q+a(\mathrm{mod}2),\; q*a=q\cdot a(\mathrm{mod}2)$.}
\end{examples}

We might refer to operations $\circ$ and $*$ as machine {\em transition} and {\em output functions}, respectively. 
From now on we would use notation $\langle Q,A,B\rangle$ to denote Mealy machine without specifying operations $\circ$ and $*$.   
Similarly we would use notation $\langle Q,A,B; q_0\rangle$ to denote initial Mealy machine. Henceforth we would use terms machine and initial machine to refer to, respectively, Mealy machine and initial Mealy machine.    

Suppose that we are given two initial machines \\
$V=\langle Q,A,B;q_0,\circ,*\rangle$ and $V'=\langle Q',A',B';q'_0,\acute\circ,\acute*\rangle$. Schematically it is shown in \ref{z56a}a.\,fig.  

\begin{figure}[h]
\scalebox{.95}{\input{aut28.pic}}
\caption{Serial composition.}
\label{z56a}
\end{figure}

We want to connect output of machine $V$ to input of machine $V'$  (shown in \ref{z56a}b.\,fig.). Clearly, in this situation, we have $v=v'$.

Suppose that $B\subseteq A'$, then for machine $V'$ input we always can use word $v=q_0*u$. Therefore word $w$ is correctly defined as
\[
w\apz q'_0\acute*(q_0*u).
\]
Our goal is to create initial machine 
$V''=\langle Q'',A,B';q''_0,\ddot\circ,\ddot*\rangle$ 
(shown \ref{z56a}c.\,fig.) such as
\[
\forall u\in A^* \; q''_0\ddot*u=q'_0\acute*(q_0*u).
\]
The new machine $V''$ is called a {\em serial} or {\em cascade composition} (or {\em connection}) of machines $V$ and $V'$ . To denote serial composition we will use notation $V\leadsto V'$.
Formally, lets define a class
\[
V\leadsto V'\apz\{\langle Q'',A,B';q''_0,\ddot\circ,\ddot*\rangle \;|\;\forall u\in A^* \; q''_0\ddot*u=q'_0\acute*(q_0*u)\}.
\]

Assume that 
\begin{itemize}
\item \quad $Q''\apz Q'\times Q$;
\item \quad $q''_0\apz(q'_0,q_0)$;
\item \quad $(q',q)\ddot\circ a\apz (q'\acute\circ q*a,q\circ a)$;
\item \quad $(q',q)\ddot*a\apz q'\acute*q*a$.
\end{itemize}

\begin{lemma} For  $\forall u\in A^*$
\begin{eqnarray*}
(q',q)\ddot\circ u&=&(q'\acute\circ q*u,q\circ u);\\
(q',q)\ddot* u&=&q'\acute* q*u.
\end{eqnarray*}
\end{lemma}

$\Box$
Proof of lemma is inductive. Let $v=ua$, where $u\in A^*$ and $a\in A$, then
\begin{eqnarray*}
(q',q)\ddot\circ v&=&(q',q)\ddot\circ ua=(q',q)\ddot\circ u\ddot\circ a=(q'\acute\circ q*u,q\circ u)\ddot\circ a\\
&=&\bigl((q'\acute\circ q*u)\acute\circ((q\circ u)*a),(q\circ u)\circ a\bigr)\\
&=&\bigl(q'\acute\circ(q*u)\acute\circ(q\circ u*a),q\circ u\circ a\bigr)\\
&=&\bigl(q'\acute\circ((q*u)\#(q\circ u*a)),q\circ ua\bigr)\\
&=&(q'\acute\circ q*ua,q\circ v)=(q'\acute\circ q*v,q\circ v);\\
(q',q)\ddot* v&=&(q',q)\ddot*ua=\bigl((q',q)\ddot* u\bigr)\#(q',q)\ddot\circ u\ddot*a\\
&=&(q'\acute*q*u)\#(q'\acute\circ q*u,q\circ u)\ddot* a\\
&=&\bigl(q'\acute*q*u\bigr)\#\bigl((q'\acute\circ q*u)\acute*((q\circ u)*a)\bigr)\\
&=&\bigl(q'\acute*(q*u)\bigr)\#\bigl((q'\acute\circ(q*u))\acute*(q\circ u*a)\bigr)\\
&=&q'\acute*\bigl((q*u)\# (q\circ u*a)\bigr)\\
&=&q'\acute*(q*ua)=q'\acute*q*v. \quad \rule{2mm}{2mm}
\end{eqnarray*}

\begin{corollary}
$V''\in V\leadsto V'$.
\end{corollary}

\section{Sequential functions}

\begin{definition} A total mapping $f:A^*\to B^*$ is called a sequential function if 
\begin{itemize}
\item[\textnormal{(i)}] $\forall u\in A^*\;|u|=|f(u)|;$
\item[\textnormal{(ii)}] $u\in \textnormal{Pref}(v)\Rightarrow f(u)\in\textnormal{Pref}(f(v))$.
\end{itemize}
\end{definition}

\begin{corollary} For all sequential functions, we have that if
\[u\in\textnormal{Pref}(v)\cap\textnormal{Pref}(w),\]
then
\[f(u)\in\textnormal{Pref}(f(v))\cap\textnormal{Pref}(f(w)).\] 
\end{corollary}

It states that if words $u$ and $v$ have matching prefixes of length $k$, then words $f(u)$ and $f(v)$ have matching prefixes of length $k$.

$\Box$ Suppose that $u\in\textnormal{Pref}(v)\cap\textnormal{Pref}(w)$, then accordingly with the definition of sequential function $f(u)\in\textnormal{Pref}(f(v))$ and $f(u)\in\textnormal{Pref}(f(w))$.~\rule{2mm}{2mm}

\begin{definition} Let $f:A^*\to B^*$ be a sequential function and $u\in A^*$, then $f_u(v)$ define a prefix of mapping $f(uv)$ with length $|v|$. The mapping $f_u$ is called a quotient of sequential function $f$.
\end{definition}

\begin{corollary} \qquad$f(uv)=f(u)f_u(v)$.
\end{corollary}

\begin{lemma}\label{l_4.2.5} \qquad $uw\in\textnormal{Pref}(uv)\Leftrightarrow w\in\textnormal{Pref}(v)$.
\end{lemma}

$\Box \Rightarrow  $ If $uw\in\textnormal{Pref}(uv)$, then there exists word $w'$ such as $uww'=uv$. Hence $ww'=v$ and therefore $w\in\textnormal{Pref}(v)$.

$\Leftarrow  $  If $w\in\textnormal{Pref}(v)$, then exists word $w'$ such as $ww'=v$. Hence $uww'=uv$ and therefore $uw\in\textnormal{Pref}(uv)$.~\rule{2mm}{2mm}

\begin{proposition} The quotient $f_u$ is a sequential function.
\end{proposition}

$\Box$ Accordingly with the definition of quotient $|v|=|f_u(v)|$.  In the same time, if $w\in\textnormal{Pref}(v)$, then [by lemma\ref{l_4.2.5}.] $uw\in\textnormal{Pref}(uv)$, therefore 
\[f(u)f_u(w)=f(uw)\in\textnormal{Pref}(f(uv)).\]
It holds that $\textnormal{Pref}(f(uv))=\textnormal{Pref}(f(u)f_u(v))$, therefore $f(u)f_u(w)\in\textnormal{Pref}(f(u)f_u(v))$. Hence [\ref{l_4.2.5}. Lemma] 
 $f_u(w)\in\textnormal{Pref}(f_u(v))$.~\rule{2mm}{2mm}

\begin{lemma}
If mapping $f:A^*\to B^*$ is a sequential function, then
\[
f_u(vw)=f_u(v)f_{uv}(w).
\] 
\end{lemma}

$\Box$  Suppose that $w\in A^*$, then we have
\begin{eqnarray*}
f(uvw)&=&f(u)f_u(vw),\\
f(uvw)&=&f(uv)f_{uv}(w)=f(u)f_u(v)f_{uv}(w).
\end{eqnarray*}
Hence $f_u(vw)=f_u(v)f_{uv}(w).~\rule{2mm}{2mm}$

\begin{proposition}
Let mapping $f:A^*\to B^*$ be a sequential function. If $f_u=f_{u'}$, then 
\[
\forall v\in A^*\;f_{uv}=f_{u'v}.
\] 
\end{proposition}

$\Box$  Let $w\in A^*$,  then 
\[f_u(v)f_{uv}(w)=f_u(vw)=f_{u'}(vw)=f_{u'}(v)f_{u'v}(w).\]
We have if $|f_u(v)|=|f_{u'}(v)|$, then $f_{uv}(w)=f_{u'v}(w)$. Hence $f_{uv}=f_{u'v}$.~\rule{2mm}{2mm}

\begin{proposition} \label{a_5.4.9}
Let mappings $f:A^* \to B^*$ and $g:B^* \to C^*$ be sequential functions, then
\[
\forall u\in A^*\; (g\circ f)_u=g_{f(u)}\circ f_u.
\]
\end{proposition}

$\Box$ 
We have 
\[
g(f(ux))=g(f(u)f_u(x))=g(f(u))\#g_{f(u)}(f_u(x)).
\]
Hence $(g\circ f)_u(x)=g_{f(u)}(f_u(x)). \; \rule{2mm}{2mm}$

\section{Restricted sequential functions}
\begin{definition} 
Let mapping $f:A^*\to B^*$ be a sequential function. The function $f$ defines set 
\[
Q_f\apz\{f_u|u\in A^*\},
\]
where $f_u$ is a quotient of $f$. The function $f$ is called a restricted sequential function (or sequential function on restricted domain) if the set $Q_f$ is finite. 
\end{definition}

\begin{theorem} For each function $f:A^*\to B^*$ there exists a initial machine $V=\langle Q_f, A, B; q_0\rangle$ such as
\[
\forall v\in A^*\; f(v)=q_0* v.
\]
\end{theorem}

We might say that output function of machine $V$ is equal with function $f$.

$\Box$ Let $Q_f=\{q_0,q_1,\ldots,q_k\}$, where element $q_0$ is a set containing quotient $f_\lambda$. 
Suppose that $f_{u}\in q$, $a\in A$, $f_{ua}\in q'$ and $f_{u}(a)=b$, then
\[
q\circ a\apz q',\qquad q*a\apz b.
\] 

(i) \qquad $f_u\in q\Rightarrow f_{uv}\in q\circ v.$

$\triangledown$ If $v=\lambda$, then $f_{u\lambda}=f_u\in q=q\circ\lambda$.

Further proof is inductive, given that $f_{uv}\in q\circ v$. Hence accordingly with the definition of the transition function $\circ$ of machine
\[
f_{uva}\in (q\circ v)\circ a=q\circ va.\quad~\vartriangle
\]

(ii) \qquad $f_u\in q\Rightarrow q*v=f_u(v).$

$\triangledown$ If $v=\lambda$, then $q*\lambda=\lambda=f_u(\lambda)$.

Further proof is inductive, given that $q*v=f_u(v)$. Since $f_{uv}\in q\circ v$, then accordingly with the definition of output function $*$ of machine $f_{uv}(a)=q\circ v*a$.
Hence 
\[q*va=q*v \#q\circ v*a=f_u(v)f_{uv}(a)=f_u(va).\quad~\vartriangle\]

We note that $f_\lambda=f$ hence $f(v)=f(\lambda)f_\lambda(v)=f_\lambda(v)$. Accordingly with the definition of element $q_0$ we have $f_\lambda\in q_0$ hence
\[f(v)=f_\lambda(v)=q_0*v.\quad~\rule{2mm}{2mm} \]

Let $P^A$ be a set, elements of which are all possible restricted sequential functions $f: A^* \to A^*$. Serial composition of two Mealy machines is a Mealy machine, 
thereby composition $f(g(u))$ of two restricted sequential functions 
\[
f: A^* \to A^*, \quad g: A^* \to A^*
\]
is restricted sequential function. Thereby we have proved that $P^A$ is a semigroup. The operation of semigroup is composition of restricted sequential functions.

In particular, when $A=\overline{0,1}^*$   we would use simpler notation  $P^1$ for denoting semigroup $P^{\overline{0,1}}$. 

\section{Group AS$_2$}

Let $\langle G,\circ\rangle$ be a monoid. An element $x \in G$ {\em has inverse} (is {\em invertible}) if 
\[
\exists y\in G\;(y\circ x=\lambda=x\circ y);
\]
where $\lambda$ is neutral element of monoid $G$. Element $y$ is called {\em dual} (or {\em inverse}) of element $x$ and is usually denoted as $x^{-1}$.

\begin{definition} A monoid $G$, each element of whom is invertible is called a group.
\end{definition}

\begin{proposition}
The set 
\[
AS_2\apz \{f\in P^1\,|\, f \textnormal{ is bijection}\}
\]
is a group in which the group operation is composition of restricted sequential functions.
\end{proposition}

$\Box$
Let $f\in AS_2$, then there exists a Mealy machine  
\[
V=\langle Q, \overline{0,1}, \overline{0,1};q_0, \circ, *\rangle
\] 
and for $\forall u\in\overline{0,1}^*\;f(u)=q_0*u$.

We henceforth use fallowing notation: $\bar 0\apz 1$ and $\bar 1\apz 0$.
Let $V'=\langle Q, \overline{0,1}, \overline{0,1};q_0, \acute \circ, \acute *\rangle$ be a new Mealy machine, where 
\begin{eqnarray*}
q\acute*a&\apz & \mbox{} \quad q*a,\\
q\acute\circ a&\apz& \begin{cases}
q\circ a,& {\rm{if}} \; \; q*a=a;\\
q\circ \bar a, & {\rm{if}} \;\; q*a=\bar a.
\end{cases}
\end{eqnarray*}
The Mealy machine $V'$ define restricted sequential function $g(u)\apz q_0\acute* u$.

We need to prove that $g=f^{-1}$, i.e., that $g$ is inverse function of $f$.

(i) Suppose that $Q'\apz\{q\,|\,\exists u\in\overline{0,1}^*\;q=q_0\circ u\}$. Lets prove that
\begin{eqnarray}\label{f5.1}
\forall q\in Q'\forall a\in\overline{0,1}\; q*(q*a)=a.
\end{eqnarray}

Let $q=q_0\circ u$ and thus
\begin{eqnarray*}
f(ua)&=&q_0*ua=q_0*u\#q_0\circ u*a=q_0*u\#q*a,\\
f(u\bar a)&=&q_0*u\bar a=q_0*u\#q_0\circ u*\bar a=q_0*u\#q*\bar a.
\end{eqnarray*}
Given that $f$ is a bijection, we have that
\begin{eqnarray}\label{f5.2}
q*a&\ne&q*\bar a.
\end{eqnarray}
\begin{itemize}
\item If $q*a=a$, then $q*(q*a)=q*a=a$.
\item If $q*a=\bar a$, then by (\ref{f5.2}) we have that $q*\bar a= a$. Hence
$q*(q*a)=q*\bar a=a$.
\end{itemize}
This proves that $fg(a)=a$.

(ii) Suppose that $q\in Q'$ and $a\in\overline{0,1}$.
\begin{itemize}
\item If $q*a=a$, then $q\acute\circ a=q\circ a$. Hence
\[
q\acute\circ(q*a)=q\acute\circ a=q\circ a.
\]
\item If $q*a=\bar a$, then by (\ref{f5.2}) $q * \bar a=a$. Hence
\[
q\acute\circ(q*a)=q\acute\circ\bar a=q\circ(q*\bar a)=q\circ a.
\]
\end{itemize}

This proves that
\begin{eqnarray}\label{f5.3}
q\acute\circ (q*a)&=&q\circ a.
\end{eqnarray}

(iii) Suppose that $q=q_0\circ u$ and 
\[
f_q: \overline{0,1}^* \to \overline{0,1}^*: v\mapsto q*v,
\]
then $f(uv)=q_0*uv=q_0*u\#q_0\circ u*v=f(u)f_q(v)$. Suppose that $w\in\overline{0,1}^*$ and we have that $f$ is a surjection therefore there exists $u_0$ such as
\[
f(u_0)=f(u)w. 
\]
Hence $u_0=u_1u_2$, where $|u_1|=|f(u)|$ and $|u_2|=|w|$. Then there exists $u_2$ such as $f(u_1u_2)=f(u)w$.
Then $u_1=u$ because $f$ is an injection. This proves $f(u)w=f(uu_2)=f(u)f_q(u_2)$. Hence $w=f_q(u_2)$, and $f_q$ is a surjection.

If $f_q(u_3)=w$, then $f(uu_3)=f(u)w$. Function $f$ is an injection, thereby $u_3=u_2$. This proves that function $f_q$ is an injection.

Finally, we conclude that $f_q$ is a bijection for each $q\in Q'$.

(iv) Suppose that $f'_q$ is a restricted sequential function defined by Mealy machine $V'$ such as
\[
f'_q: \overline{0,1}^* \to \overline{0,1}^*: u\mapsto q\acute * u.
\]
In section (i) we proved (equation (\ref{f5.1})), i.e., $f_qf'_q(a)=a$, if $a\in\overline{0,1}$.
Further proof is inductive, given that for all $q\in Q'$ and all words $u$ of length $n$ holds $f_qf'_q(u)=u$.

Suppose that $a\in\overline{0,1}$, then 
\begin{eqnarray*}
f_qf'_q(au)&=&q*q\acute* au = q*(q\acute* a\# q\acute\circ a\acute* u) =q*(q*a\#q\acute\circ a\acute* u)\\
&=&q*q*a\#q\circ(q*a)*(q\acute\circ a\acute* u)\\
&=&a\#q\acute\circ a*(q\acute\circ a\acute* u)=a\#f_{q\acute\circ a}f'_{q\acute\circ a}(u)=au.
\end{eqnarray*}
We note that $q\acute\circ a=q\circ(q*a)\in Q'$. This concludes inductive part of the proof.
We note that $f(u)=f_{q_0}(u)$ and $f^{-1}(u)=f'_{q_0}(u)$. 

(v) Let $f$ and $g$ be restricted sequential functions from the set $AS_2$, Then there exists Mealy machines
\[
V=\langle Q, \overline{0,1}, \overline{0,1};q_0, \circ, *\rangle,\;
V'=\langle Q', \overline{0,1}, \overline{0,1};q'_0, \acute\circ, \acute*\rangle
\]
such that for $\forall u\in\overline{0,1}^*$
\[
f(u)=q_0*u\;\wedge\;g(u)=q'_0\acute* u \,
\]
The serial composition $ V\leadsto V'$ is machine, which realizes composition of functions $f$ and $g$. In other words, if   
$\breve V =\langle \breve Q, \overline{0,1}, \overline{0,1}; \breve q_0, \breve\circ, \breve*\rangle\in  V\leadsto V' $, then 
\[
\forall u\in\overline{0,1}^*\; gf(u)=\breve q_0\breve* u.
\]

\begin{figure}[h] \label{aut35}
\input{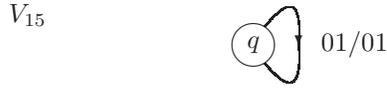}
\caption{Implementation of identity function.}
\label{z54}
\end{figure}

(vi) The identity function $\mathbb{I}: \overline{0,1}^* \to \overline{0,1}^* : u\mapsto u$ is a restricted sequential function. The machine implementing identity function given as $V_{15}$ (shown in \ref{z54}fig.).

(vii) We note that compositions of functions is associative.  Finally, we have proven (in (iv) -- (vi)) that $AS_2$ is a group.
\rule{2mm}{2mm}

Suppose that $ X\subseteq G$, then 
\[\langle X\rangle\apz\bigcap_{X\subseteq H\le G}H,\]
i.e., we consider intersection of all subgroups of group $G$ containing set $X$.

{\bf Notification}.
In this situation similarly as in case of the semigroups we use notation $\langle X\rangle$.
\begin{definition}
The group $\langle X\rangle$ is called subgroup generated by set $X$.
\end{definition}

If $\langle X\rangle=G$, then set $X$ is called the{\em set of generators} of group $G$ . Elements of set $X$ are called {\em generators}. 
If $X=\{x_1,x_2,\ldots,x_n\}$, then the following notification is used
 \[\langle x_1,x_2,\ldots,x_n\rangle\apz\langle X\rangle.\] 
In this case (set $X$ --- finite), group $G$ is {\em finitely generated.}

\begin{definition}\label{d3.8.4}
A group $G$ is called cyclic group if
\[\exists g\in G\;\forall x\in G\;\exists n\in\mathbb{Z}\quad x=g^n.\]
\end{definition}

In this case an element $g$ is called a {\em generator} of cyclic group $G$ .

\begin{definition}
Let $\langle G,\odot\rangle$ be group. The cardinality of set $G$, denoted as $|G|$, is called the order of group $G$.
\end{definition}

A group $G$ is called a {\em finite group} if its cardinality $|G|$ is a natural number, denoted $|G|<\aleph_0$. In opposite case, a group $G$ is called an {\em infinite group} and denotated as $|G|\ge\aleph_0$. 

\begin{definition}
Let $a$ be an element of group $G$, then subgroups $\langle a\rangle$ order is called the order of element $a$.
\end{definition}

We use notification $o(a)$ to denote order of element, i.e., $o(a)\apz|\langle a\rangle|$.

\section{Elements of graph theory}

\begin{definition}
A 2-sorted algebra $\langle V,E,s,t \rangle$ is called pseudograph if:\\ 
\textnormal{(i)} $V$ is nonempty set;\\
\textnormal{(ii)} $s,t$ are total mappings $E\overset{s}{\rightarrow}V,E\overset{t}{\rightarrow}V$, respectively.
\end{definition}
The elements of the set $V$ are called {\em vertices},the elements of the set $E$ are called {\em arcs}.
If sets $V$ and $E$ are finite then the pseudograph $G$ is called finite. In this case to denote a directed arc $l$ we use notation $(s(l),t(l))$. The mapping $s$ is called an {\em initial vertex}, while the mapping $t$ is called a {\em terminal vertex}. We say that an arc $l$ is {\em incident out} of the vertex $s(l)$ and {\em incident in} to the vertex $t(l)$. 
\begin{definition}
A pseudograph $G^{'}=\langle V^{'},E^{'},s^{'},s^{'}\rangle$ is called a sub-pseudograph of pseudograph $G=\langle V,E,s,t\rangle$ if :\\
\textnormal{(i)} $V^{'}\subseteq V$; \\
\textnormal{(ii)} $E^{'}\subseteq \left \{ l \in E  |  s(l) \in V^{'}\wedge t(l) \in V^{'} \right \}$; \\
\textnormal{(iii)} $s^{'}=s|_{E^{'}} \wedge t^{'}=t|_{E^{'}}$.
\end{definition}
\begin{definition}
a pseudograph $\langle V,E,s,t \rangle$ is called an oriented graph, if :\\
\textnormal{(i)} $\forall l\in E \quad s(l)\ne t(l)$\\
\textnormal{(ii)} $\forall l_1l_2\in E  \quad (s(l_1)=s(l_2)\wedge t(l_1)=t(l_2)\Rightarrow l_1=l_2).$
\end{definition}
This means that oriented graph can be given by two sets $V$, $E\subseteq V^2$. Usually in this case notation $G(V,E)$ is used.

\begin{definition}
An oriented graph is called a graph (non oriented graph), if
$\forall \,l_1\in E,\exists\, l_2\in E: (s(l_1)=t(l_2)\wedge s(l_2)=t(l_1)).$
In this case the pair of symmetrical arcs $l_1,l_2$ is identified by the set $\{s(l_1),s(l_2)\}$. This pair is called an edge connecting vertexes $s(l_1)$ and $s(l_2)$.
\end{definition}
To denote edge $\{s(l_1),s(l_2)\}$ we use notification $[s(l_1),s(l_2)]$.

\begin{definition}
A $2n+1$ tuple $c=(v_{0},l_{0},v_{1},l_{1},...,l_{n},v_{n+1})$ is called a walk in pseudograph $\langle V,E,s,t \rangle$ if\\
\textnormal{(i)} $\forall i \in \overline{0,n+1}, v_{i}\in V$;\\
\textnormal{(ii)} $\forall j \in \overline{0,n}, l_{j}\in E$;\\
\textnormal{(iii)} $\forall j \in \overline{0,n}, s(l_{j})=v_{j} \wedge t(l_{j})=v_{j+1}$. 
\end{definition}

We say that a walk $c=(v_{0},l_{0},v_{1},l_{1},...,l_{n},v_{n+1})$ of {\em length} starts at the vertex $v_{0}$ and ends at the vertex $v_{n+1}$, connecting vertex $v_{i}$ with $v_{j}$, $i<j, i,j \in \overline{0,n+1}$.
From the definition of graph, for each walk $c=(v_{0},l_{0},v_{1},l_{1},...,l_{n},v_{n+1})$ there exists a {\em reverse} walk
\[
c^{'}=(v_{n+1},l_{n},v_{n},l_{n-1},...,l_{0},v_{0}),
\]
therefore we may say that the walks $c$ connects vertices $v_{i}$ and $v_{j}$, $i,j \in \overline{0,n+1}$.

Let $m=(m_1,m_2,\ldots, m_k)$ be a walk, here the integer $k$ is called {\em extended length} of walk $m$. A vertex $v$ occurs in a walk $m$, if there $\exists j \; v=m_j$. An arc $l$ occurs in a walk $m$ if $\exists i \; l=m_i$. An edge $e=\{s(l_1),s(l_2)\}$ occurs in a walk if $\exists i \; l_{1}=m_i \vee l_{2}=m_i$. The walk $m^{'}=(m^{'}_1,m^{'}_2,\ldots, m^{'}_l)$ is a {\em part} of walk $m=(m_1,m_2,\ldots, m_k)$ if $\exists j \;l\leq k-j \wedge \forall i \in \overline{1,l} m^{'}_{i}=m_{j+i}$.

\begin{definition}
A walk is called 
\begin{itemize}
\item closed, if it starts and ends at the same vertex;
\item trail, if it is without repeated arcs (and edges);
\item cycle, if it is nonempty closed trail;
\item simple trail (or a path), if it is without repeated vertices.
\end{itemize}
\end{definition}

\begin{proposition}\label{a5.8.7}
If there exists a walk connecting vertices $v_1$ and $v_2$, then there exists a trail connecting $v_1$ and $v_2$.
\end{proposition}

$\Box$ The proof is inductive by induction on extended length $k$ of walk $m=(m_1,m_2,\ldots, m_k)$. If the extended length of walk is 1 then it is a trail.

If $m$ is not a trail, then there exists an edge $l$ such that $l=m_i=m_{i+j}$, where $1<i<i+j<k$.
Then for vertices $m_{i-1}=m_{i+j-1}$, $m_{i+1}=m_{i+j+1}$. 
The extended length of walks 
\[
(m_1,\ldots, m_{i-1}, l, m_{i+j+1}, m_{i+j+2},\ldots, m_k)
\]
is less than $k$, therefore accordingly to the inductions hypothesis there exist a trail connecting $v_1$ and $v_2$.
\rule{2mm}{2mm}

\begin{proposition}\label{a5.8.7a}
If there exists two different trails connecting vertices $v_1$ and $v_2$, then there exist a cycle.
\end{proposition}

$\Box$
Suppose that $(v_1, \alpha_1,\alpha_2 ,\ldots, \alpha_k, v_2)$ and $(v_1, \beta_1,\beta_2 ,\ldots, \beta_n, v_2)$ are two different trails connecting vertices $v_1$ and  $v_2$.
Suppose that
\[
\alpha_1=\beta_1, \alpha_2=\beta_2 ,\ldots, \alpha_s=\beta_s, \alpha_{s+1}\ne \beta_{s+1},
\]
then $\alpha_s$ is a vertex and 
\[
(\alpha_s,\alpha_{s+1},\ldots, \alpha_k, v_2), \quad (\beta_s, \beta_{s+1},\ldots, \beta_n, v_2)
\]
are two different trails connecting vertex $\alpha_s=\beta_s$ and vertex $v_2$.

Suppose that $\alpha_t$ is the first element of sequence
\[
\alpha_{s+1}, \alpha_{s+2}, \ldots, \alpha_k,
\]
matching with some element $\beta_\tau$ in  $\{\beta_{s+1},\beta_{s+2},\ldots, \beta_n\}$. Then
\[
(\alpha_s,\alpha_{s+1},\ldots, \alpha_t=\beta_\tau, \beta_{\tau-1},\ldots,\beta_{s+1},\beta_s=\alpha_s)
\]
is a cycle, otherwise
\[
(\alpha_s,\alpha_{s+1},\ldots,\alpha_k,v_2,\beta_n,\ldots,\beta_{s+1},\beta_s)
\]
is a cycle.
\rule{2mm}{2mm}

\begin{definition}
An edge $[v_1,v_2]$ of graph is called to be incident to vertex $v_1$, and also to vertex $v_2$.  Such vertices $v_1$ and $v_2$ are called neighboring vertices. Thus vertex $v_1$ is a neighbor of vertex $v_2$ and vice versa.
The number of edges incident to a vertex $v$ is called the degree of vertex.
\end{definition}
In this case we use notification ${\rm deg}\,(v)$. A vertex $v$ is called {\em isolated vertex}, if ${\rm deg}\,(v)=0$. A vertex $v$ is called a {\em leaf} if ${\rm deg}\,(v)=1$.

\begin{proposition}\label{a5.8.9}
Let $G(V,E)$ be a finite graph, then 
\[
\sum_{v\in V}{\rm deg}\,(v)=2|E|
\]
\end{proposition}

$\Box$ As each edge is incident to two vertices, the sum of all degrees of vertices is twice the number of edges $2|E|$.
\rule{2mm}{2mm}

\begin{corollary}\label{s5.8.10}
If $G(V,E)$ is a graph in which the degree of each vertex $v$ is ${\rm deg}\,(v)\ge 2$, then for sets $V,E$ holds $|E|\ge|V|$.
\end{corollary}

$\Box$
$2|E|=\sum\limits_{v\in V}{\rm deg}\,(v)\ge 2|V|$. Hence $|E|\ge|V|$.
\rule{2mm}{2mm}

\begin{definition}
A graph is called connected if for each two vertices $v_1, v_2$ there exist a walk $m$ starting at $v_1$ and ending at $v_2$.
Otherwise graph is called disconnected. 
An edge $[a,b]$ of connected graph G(V,E) is called a bridge if graph $G(V,E\setminus \{[a,b]\})$ is disconnected.
A connected graph $G(V_1,E_1)$ is called {\em connected component} of $G(V,E)$ if it is maximal connected subgraph.
\end{definition}

\begin{lemma}\label{l5.8.12}
If $|E|=|V|-2$, then the graph $G(V,E)$ is disconnected.
\end{lemma}

$\Box$ Proof done by induction on $n=|V|$. If $n=2$, then $n-2=0$ hence graph with two vertices and no edges is disconnected.

Let $G(V,E)$ be a graph with $n+1$ vertices and $n-1$ edges. By (corollary \ref{s5.8.10}) there exist vertex $v$ with degree 1 or 0.
If vertex $v$ is isolated then graph $G(V,E)$ is disconnected. Otherwise, if $v$ is a leaf, we exclude vertex $v$ and edge incident to it from the graph. The obtained graph $G'$ has $n$ vertices and $n-2$ edges. From the induction hypothesis graph $G'$ is not connected. Therefore there exist two vertices $v_1$ and $v_2$ for whom there is no trail connecting them. 
If such trail $\varkappa$ would exists in the graph $G(V,E)$ then vertex $v$ would occur in it. We have that $v_1\ne v\ne v_2$ hence walk $\varkappa$ can't be a trail because vertex $v$ is a leaf.
Hence graph $G(V,E)$ is disconnected. 
\rule{2mm}{2mm}

\begin{definition}\label{d5.8.13a}
A graph $P$ is called an underlying graph of oriented graph $G$, if its obtained by replacing all arcs of graph $G$ with edges. 
\end{definition}
Replacement operation involves assigning to each arc in the original graph an arc symmetrical to it.

\begin{definition}\label{d5.8.13}
A connected graph without cycles is called a tree. An oriented graph is called a oriented tree, if its underlying graph is a tree.
An oriented tree is called a rooted oriented tree, if there exist a vertex $\sigma$ such that for every other vertex $v$ there exist a walk connecting $\sigma$ and $v$. Such vertex $\sigma$ is called a root.
\end{definition}

\begin{proposition} \label{a5.8.16}
The fallowing properties of graph T are equivalent:
\begin{enumerate}
\item $T=T(V,E)$ is a tree;
\item $T$ is connected and each edge is a bridge;
\item $T$ is connected and contains $n-1$ edge;
\item $T$ is without cycles and contains $n-1$ edge;
\item Each two vertices of $T$ are connected only one trail;
\item $T$ is without cycles, but connecting any two non neighboring vertices with a new edge would create a cycle.
\end{enumerate}
\end{proposition}

$\Box \; 1.\Rightarrow 2.$ The graph $T$ is connected by definition, by the exclusion of an edge we cannot maintain connected graph, otherwise, there should be a cycle.\\
$2. \Rightarrow 3.$ The proof is done by an induction on $n$. Let exclude one edge, say, $[v_1,v_2]$. We obtain a new graph $T(V,E_0)$, where $E_0\apz E\setminus\{[v_1,v_2]\}$.
Let define a subset of vertex set of the graph $T(V,E_0)$:
\begin{eqnarray*}
V_1&\apz&\{v\,|\, \textnormal{there exist a walk connecting vertex } v_1 \textnormal{ and } v \},\\
V_2&\apz&\{v\,|\, \textnormal{there exist a walk connecting vertex } v_2 \textnormal{ and } v \}.
\end{eqnarray*}

(i) The graph $T(V,E_0)$ is disconnected, therefore $V_1\cap V_2=\emptyset$. Otherwise there would exist a vertex $v\in V_1\cap V_2$, implying existence of walks $\mu_1$, $\mu_2$ connecting vertex $v_1$ to $v$ and, respectively, vertex $v$ to $v_2$. Therefore there would exist a walk $m=(m_1,m_2,\ldots,m_s)$ connecting vertices v$v_1$ and $v_2$.

Let $w_1, w_2$ be two arbitrary chosen vertices of the graph $T(V,E)$. As the graph $T$ is connected, there exist a walk connecting $w_1$ and $w_2$. Therefore (by proposition \ref{a5.8.7}) there exist
a trail $\varkappa=(\varkappa_1,\varkappa_2,\ldots, \varkappa_\nu) $ connecting $w_1$ and $w_2$. If the edge $[v_1,v_2]$ is not occurring in the trail $\varkappa$, then $\varkappa$ is a trail in graph
$T(V,E_0)$ connecting vertices $w_1$ and $w_2$. If trail $\varkappa$ contains edge $[v_1,v_2]$, then exists index $i$ such as $[v_1,v_2]=\varkappa_i$. Hence for the part of trail we have
\begin{eqnarray*}
(\varkappa_{i-1},\varkappa_i,\varkappa_{i+1})&=&(v_1,[v_1,v_2],v_2),\\
{\rm or}\hspace{40mm} \\
(\varkappa_{i-1},\varkappa_i,\varkappa_{i+1})&=&(v_2,[v_1,v_2],v_1).
\end{eqnarray*}
In first case 
\[
(\varkappa_1,\ldots,\varkappa_{i-1},m_2, m_3,\ldots,m_{s-1},\varkappa_{i+1},\varkappa_{i+2},\ldots,\varkappa_\nu)
\]
is a walk in graph $T(V,E_0)$ connecting vertices $w_1$ and $w_2$. In second case 
\[
(\varkappa_1,\ldots,\varkappa_{i-1},m_{s-1}, m_{s-2},\ldots,m_2,\varkappa_{i+1},\varkappa_{i+2},\ldots,\varkappa_\nu)
\]
is walk in graph $T(V,E_0)$ connecting vertices $w_1$ and $w_2$. Therefore graph $T(V,E_0)$ is connected. Thus we have arrived at a contradiction!

(ii) Lets prove that $V_1\cup V_2=V$. Suppose that $v\in V_1\cup V_2$, thee there exist a trail
\[
\alpha = (\alpha_1,\alpha_2,\ldots, \alpha_t),
\]
in graph $T$ connecting vertices $v$ and $v_2$. If edge $[v_1,v_2]$ occurs in trail $\alpha$, then $v\in V_2$. If edge $[v_1,v_2]$ occurs in trail $\alpha$, then there exist index $i$ such that $\alpha_i=[v_1,v_2]$.
Hence for the part of trail
\begin{eqnarray*}
(\alpha_{i-1},\alpha_i,\alpha_{i+1})&=&(v_1,[v_1,v_2],v_2),\\
{\rm or}\hspace{40mm} \\
(\alpha_{i-1},\alpha_i,\alpha_{i+1})&=&(v_2,[v_1,v_2],v_1).
\end{eqnarray*}
In first case $(\alpha_1,\alpha_2,\ldots,\alpha_{i-1})$ is a trail in graph  $T(V,E_0)$ connecting $v$ and $v_1$, therefore $v\in V_1$. In second case  $(\alpha_1,\alpha_2,\ldots,\alpha_{i-1})$
is a trail in $T(V,E_0)$ connecting $v$ and $v_2$, therefore $v\in V_2$.

Hence graph $T(V,E_0)$ contains two connected components\\
$T_1(V_1,E_1)$ and $T_2(V_2, E_2)$, $E{_{1}}=E_{0}|_{V_{1}\times V_{1}}$, $E{_{2}}=E_{0}|_{V_{2}\times V_{2}}$. As each edge in graph $T$ is a bridge, then  in graphs $T(V_1,E_1)$ and $T(V_2,E_2)$ all edges are bridges.

According to the induction hypothesis: $|E_1|=|V_1|-1$ and $|E_2|=|V_2|-1$. 
As a result we have
\[
|E|=|E_1|+|E_2|+1=|V_1|-1+|V_2|-1+1= |V|-1.
\]

$3. \Rightarrow 4.$ 
If graph would contain a cycle then excluding some edge belonging to this cycle would not make graph disconnected. The graph $T_0$ obtained in this way would be connected and would contain $n-2$ edges and $n$ vertices. This is in contradiction by lemma \ref{l5.8.12}.

$4. \Rightarrow 5.$ From given is clear that no two vertices are connected with more than one trail, otherwise there would exist a cycle.
We need to prove that graph is connected. Suppose that graph is disconnected, and contains multiple connected components. 
Each connected component is without a cycles, therefore (by definition \ref{d5.8.13}) they are trees. Hence 3. condition is true, i.e., each connected component $G(V_i,E_i)$ has property $|E_i|=|V_i|-1$. Hence, if  number of connected component are $k$:
\begin{eqnarray*}
&&\; \; \qquad \sum_{i=1}^k |E_i|=\sum_{i=1}^k|V_i|-k,\\
n-1&=&|E|=\sum_{i=1}^k |E_i|=\sum_{i=1}^k|V_i|-k=|V|-k=n-k.
\end{eqnarray*}
Therefore $k=1$, i.e., we have only one connected component.

$5. \Rightarrow 6.$ If there exist a cycle in graph $T$, then at least two distinct vertices could be connected with two different trails. This leads to contradiction and therefore graph $T$ is without cycles. 

Suppose that vertices $v_1,v_2$ aren't neighbors, then exists a trail 
\[
(v_2,\varkappa_1,\ldots,\varkappa_k,v_1),
\]
connecting $v_2$ and $v_1$. Let add a edge $[v_1,v_2]$, then 
\[
(v_2,\varkappa_1,\ldots,\varkappa_k,v_1,[v_1,v_2],v_2)
\]
is a cycle.

$6. \Rightarrow 1.$ Suppose that $v_1,v_2$ aren't neighbors, but adding edge $[v_1,v_2]$ create a cycle. Therefore $[v_1,v_2]$ occurs in a cycle.
Let suppose that this new cycle is 
\[
(v_1,[v_1,v_2],v_2,\varkappa_1,\ldots,\varkappa_k,v_1),
\] 
the $(v_2,\varkappa_1,\ldots,\varkappa_k,v_1)$ is a walk connecting $v_2$ and $v_1$. Therefore $T$ is a connected graph, additionally $T$ is without cycles. 
Hence, accordingly to the definition of tree, $T$ is a tree.
\rule{2mm}{2mm}

\begin{definition} The integer ${\rm deg}^-(v)\apz |\{l|t(l)=v\}|$ is called the indegree of vertex $v$ of pseudograph $\langle V, E, s,t\rangle$.
The integer ${\rm deg}^+(v)\apz |\{l|s(l)=v\}|$ is called the outdegree of vertex $v$ of pseudograph $\langle V, E, s,t\rangle$.
\end{definition}
A pseudograph is called {\em out--$p$--regular}, if outdegree of all vertices is $p$.
Pseudograph $G$ is {\em out--regular}, if there exists such integer $p$, that $G$ is out--$p$--regular pseudograph. A rooted oriented tree is called regular ($p$--regular), if it is out--$p$--regular.

\begin{corollary}\label{s5.8.17}
$\sum\limits_{v\in V}{\rm deg}^+(v)=|E|$.
\end{corollary}

$\Box$
Each arch starts with a vertex.
\rule{2mm}{2mm}

\begin{definition}
A graph $G(V,E)$ is called vertex infinite, if its vertex set $V$ is infinite.
\end{definition}

\begin{corollary}
If $p>0$, then each out--$p$--regular rooted tree is infinite.
\end{corollary}
$\Box$
If infinite rooted tree $G(V,E)$ is out--$p$--regular, then
\[
p|V|=\sum_{v\in V}{\rm deg}^+(v)\overset{S\ref{s5.8.17}}{=}|E|.
\]
If $p>0$, then $|V|\le p|V|=|E|$, but for all finite trees holds $|E|=|V|-1$. This is a contradiction! 
\rule{2mm}{2mm}

\begin{corollary}
The indegree for root of an out--$p$--regular rooted tree is ${\rm deg}^-(\sigma)=0$.  
\end{corollary}

$\Box$ Suppose that $\sigma$ is the root of a out--$p$--regular tree and ${\rm deg}^-(\sigma)>0$, then there exist a vertex $v$ and an arch $l$ such as $s(l)=v$ and $t(v)=\sigma$.
As $\sigma $ is the root, then there exist a walk
\[
(\sigma, c_1, c_2,\ldots, c_n,v),
\]
connecting $\sigma$ with $v$. Hence
\[
(\sigma, c_1, c_2,\ldots, c_n,v,l,\sigma)
\]
is a cycle. A contradiction!
\rule{2mm}{2mm}

\begin{proposition}\label{a5.8.23}
If there exist a trail connecting vertices $v_1$ and $v_2$, then there exist a simple trail (path) connecting $v_1$ and $v_2$. 
\end{proposition}

$\Box$ Suppose that $c=(c_1,c_2,\ldots, c_n)$ is a trail connecting $v_1$ and $v_2$, i.e., $c_1=v_1$ and $c_n=v_2$. 
If $v_1=v_2$, then trail $c=(v_1)$ is simple. 

Let $v_1\ne v_2$. The proof is done by induction on length of trail $c$.
If $c$ length is 3, then it is a simple trail.
Suppose that $n>3$ and $c$ is not a simple trail, then there exist vertex $c_i$ such that $c_i=c_{i+j}$, where
\[
1\le i<i+j\le  n.
\]
We can not have both $1=i$ and $i+j=n$ in the some time because $c_1=v_1\ne v_2=c_n$. Trail
\[
(c_1,c_2,\ldots, c_i, c_{i+j+1},\ldots, c_{n-1}, c_n)
\]
connects $v_1$ and $v_2$. Length of the trail is less than n, therefore by the induction hypothesis exists a simple trail connecting $v_1$ and $v_2$.
If $i+j=n$, then trail $(c_1,c_2,\ldots, c_i)$ connects $v_1$ and $v_2$. The length of this trail also is less than n, therefore by the induction hypothesis exists simple trail connecting $v_1$ and $v_2$.
\rule{2mm}{2mm}

\begin{corollary}
If $v$ is a vertex of a out--$p$--regular rooted tree and $v$ is not a root then ${\rm deg}^-(v)=1$.  
\end{corollary}

$\Box$ Let $G$ be a $p$--regular tree and $P$ an underlying graph of $G$. Suppose that ${\rm deg}^-(v)>1$, then exists at least 2 vertices $v_1$ and $v_2$ such that $(v_1,v)$ and $(v_2,v)$ are arcs.
Suppose that $\sigma$ is the root, then exists simple trails
\[
(\sigma,\alpha_1, \alpha_2, \ldots, \alpha_k, v_1) \quad {\rm and} \quad (\sigma, \beta_1, \beta_2 ,\ldots, \beta_n, v_2).
\]

(i) Is the walk $c_1\apz(\sigma,\alpha_1, \alpha_2, \ldots, \alpha_k, v_1, (v_1,v),v)$ a simple trail?

If we assume contrary, then there exists a vertex $\alpha_i$ such as $\alpha_i=v$. Hence, if orientation of arcs is removed),
\[
(\alpha_i,\alpha_{i+1},\ldots, \alpha_k, [v_1,v],v)
\]
is a cycle in underlying graph $P$. This is a contradiction, since $P$ is a tree. Therefore $c_1$ is a simple trail.
Similarly is provable that $c_2\apz(\sigma,\beta_1, \beta_2, \ldots, \beta_n, v_2, (v_2,v),v)$  is a simple trail.

(ii) If the orientation of arcs is removed, then $c_1$ and $c_2$ are 2 distinct trails in the underlying graph $P$ connecting vertices $\sigma$ and $v$. This is a contradiction, since $P$ is a tree.
\rule{2mm}{2mm}

\begin{proposition}\label{a5.8.25}
For each vertex $v$ of a $p$--regular tree there exist only one simple trail connecting root $\sigma$ with vertex $v$.
\end{proposition}

$\Box$ Let $c_v$ be the shortest trail connecting root $\sigma$ with vertex $v$, i.e., if $c$ is a trail connecting $\sigma$ with $v$, then trail $c_v$ is not longer than $c$. As ${\rm deg^-}(\sigma)=0$, then only trail connecting $\sigma$ with $\sigma$ is of length $0$.

Further proof is inductive, given that for each vertex $v$ such as $c_v\le l$ exists only one simple trail connecting root $\sigma$ with $v$.
Suppose that $w$ is vertex such as trail $c_w=(\sigma,\alpha_1,\alpha_2,\ldots,\alpha_n, \alpha_{n+1},w )$ is of length $l+1$. \\
Here $c\apz (\sigma,\alpha_1,\alpha_2,\ldots, \alpha_n)$ is a trail of length $l$, $\alpha_{n+1}=(\alpha_n,w)$ is an arc and $\alpha_n$ is a vertex.
As ${\rm deg^-}(w)=1$, then for each trail $(\sigma,\beta_1,\ldots, \beta_k,\beta_{k+1},w)$ connecting $\sigma$ with $w$, arch $\beta_{k+1}=(\alpha_n, w)$ and vertex $\beta_k=\alpha_n$.

Therefore $(\sigma,\beta_1,\ldots, \beta_k)$ is trail connecting $\sigma$ with $\alpha_n$, i.e., 
\[
(\sigma,\beta_1,\ldots, \beta_{k-1},\beta_k)=(\sigma,\beta_1,\ldots, \beta_{k-1},\alpha_n).
\]
The length of trail $c$ is $l$, therefore length of $c_{\alpha_n}$ is not larger than $l$. From induction hypothesis, there exist only one trail connecting $\sigma$ with $\alpha_n$.
Then
\begin{eqnarray*}
c_{\alpha_n}&=&c=(\sigma,\alpha_1,\alpha_2,\ldots, \alpha_n)\\
&=&(\sigma,\beta_1,\ldots, \beta_{k-1},\alpha_n)\\
&=&(\sigma,\beta_1,\ldots, \beta_{k-1},\beta_k).
\end{eqnarray*} 
Hence
\begin{eqnarray*}
(\sigma,\beta_1,\ldots,\beta_{k-1}, \beta_k,\beta_{k+1},w)&=&(\sigma,\alpha_1,\alpha_2,\ldots, \alpha_n,\beta_{k+1},w)\\
&=&(\sigma,\alpha_1,\alpha_2,\ldots, \alpha_n,(\alpha_n,w),w)\\
&=&(\sigma,\alpha_1,\alpha_2,\ldots, \alpha_n,\alpha_{n+1},w)\\
&=&c_w.
\end{eqnarray*} 
Finally, we conclude that unity of trail is proved.
\rule{2mm}{2mm}
\begin{definition}
The integer $l$ is called level of vertex $v$ if $l$ is length of trail connecting root $\sigma$ with vertex $v$.
\end{definition}

\begin{corollary}\label{s5.8.27}
For each vertex $v$ of a $p$--regular tree exists only one level.
\end{corollary}

Suppose that $G(V,E)$ is a $p$--regular tree with root $\sigma$. Let $l(v)$ denote length of trail connecting root $\sigma$ with vertex $v$.
\[
L(n)\apz\{v\,|\, l(v)=n\}.
\]
Elements of set $L(n)$ are called $n$--th level vertices. From collorary \ref{s5.8.27} set $L(n)$ is defined uniquely and correctly.

\begin{definition}\label{a5.8.27a}
$|L(n)|=p^n$ for each $p$--regular tree $G(V,E)$.
\end{definition}

$\Box$ As $G(V,E)$ is $p$--regular, then $|L(1)|=p$. Further proof is inductive , given that $|L(n)|=p^n$. 
Suppose that $v_1,v_2\in L(n)$ and $(v_1,u_1),(v_2,u_2)$ are arcs of tree $G(V,E)$, then $u_1\ne u_2$. Otherwise, there would exist 2 distinct trails connecting root to vertex $u_1$. This is in contradiction with proposition \ref{a5.8.25}.

Suppose $G(V,E)$ is a $p$--regular, therefore for each vertex $v$ of level $n$ the outdegree ${\rm deg}^+(v)=p$. Hence
\[
|L(n+1)|=|L(n)|p=p^np=p^{n+1}. \quad \rule{2mm}{2mm}
\]

\begin{definition}
Pair of mappings $(f_1: V_1 \to V_2, f_2: E_1 \to E_2)$ is called a homomorphism of  pseudographs 
\[
\langle V_1, E_1, s_1,t_1\rangle \; \langle V_2, E_2, s_2,t_2\rangle,
\]
if for each arch $l$ of pseudograph $\langle V_1, E_1, s_1,t_1\rangle$ holds:\\
\textnormal{(i)} $f_1(s1(l))=s_2(f_2(l))$,\\
\textnormal{(ii)} $f_1(t_1(l))=s_2(f_2(l))$.
\end{definition}
If mappings $f_1$ and $f_2$ are bijections, then homomorphism is called an {\em isomorphism}. If $\langle V_1, E_1, s_1,t_1\rangle=\langle V_2, E_2, s_2,t_2\rangle$,
then homomorphism is called an {\em endomorphism}. If additionally $f_1,f_2$ are bijections, then endomorphism is called an {\em automorphism}.

\begin{definition}
A 3-sorted algebra $\langle V,E,A,s,t,i\rangle$ is called a labeled pseudograph if :\\
\textnormal{(i)} $\langle V,E,s,t\rangle$ is a pseudograph,\\
\textnormal{(ii)} $i$ is a total mapping $E\overset{i}{\rightarrow}A$.
\end{definition}

Image $i(l)\in A$ is called a label of arc $l \in E$.

\begin{figure}[h]
\input{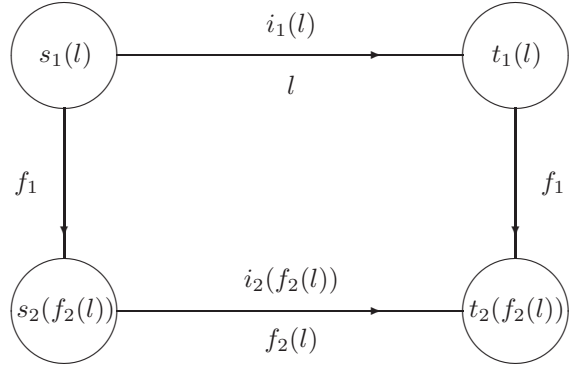}
\caption{Isomorphism of labeled pseudographs.}
\label{z69}
\end{figure}

\begin{definition}
Labeled pseudographs 
\[
\langle V_1, E_1, s_1,t_1,i_1\rangle, \; \langle V_2, E_2, s_2,t_2,i_2\rangle
\]
are isomorph labeled pseudographs, if there exist bijections
\[
f_1: V_1 \to V_2, \: f_2: E_1 \to E_2, 
\]
such that for each arch $l$ of a labeled pseudograph $\langle V_1, E_1, s_1,t_1,i_1\rangle$ holds:\\
\textnormal{(i)} $f_1(s_1(l))=s_2(f_2(l))$,\\
\textnormal{(ii)} $f_1(t_1(l))=s_2(f_2(l))$,\\
\textnormal{(iii)} $i_1(l)=i_2(f_2(l))$.
\end{definition}
Look for illustration in \ref{z69} fig.

If pseudographs $\langle V_1, E_1, s_1,t_1\rangle, \langle V_2, E_2, s_2,t_2\rangle$ are oriented graphs $G(V_1,E_1), G(V_2,V_2)$, then mapping $f_1:V_1\to V_2$ uniquely defines mapping $f_2: E_1\to E_2$ :
\[
f_2: E_1 \to E_2: (u,v)\mapsto (f_1(u), f_1(v)).
\]
Therefore for the definition of homomorphism of oriented graphs only one mapping is needed.

\begin{definition}\label{d5.8.32}
A mapping $f: V_1 \to V_2$ is called a homomorphism of (oriented) graphs 
$G(V_1,E_1), G(V_2,E_2)$, if for each arc $(u,v)\in E_1$ there exist $(f(u),f(v))\in V_2$.
\end{definition}
If $G(V_1,E_1)=G(V_2,E_2)$, then a homomorphism is called an {\em endomorphism} of oriented graph $G(V_1,E_1)$.

The situation is different with (oriented) graph isomorphism. Condition stating that $f:V_1\to V_2$ is a bijection is not enough (look \ref{z72}\,fig.).
Mapping
\[
f:u_1\mapsto v_1, u_2\mapsto v_2, u_3\mapsto v_3
\]
is a bijection of vertex sets of graphs $G_1,G_2$, therefore an homomorphism, but it does not define the bijection of sets of arcs.

\begin{figure}[h]
\input{aut49.pic}
\caption{  Nonisomorphic graphs.}
\label{z72}
\end{figure}

\begin{definition}\label{d5.8.33}
A homomorphism $f$ of (oriented) graphs $G_1$, $G_2$ is called an isomorphism if $f$ is a bijection of the vertex sets and $f^{-1}$ is a graph homomorphism of $G_2,G_1$.
\end{definition}
If $G_1=G_2$, then isomorphism is called an {\em automorphismu} of graph $G_1$. 

\begin{definition}\label{d5.8.34}
Let $G(V_1,E_1), G(V_2,E_2)$ be oriented rooted trees. Then a homomorphism $f$ of oriented graphs $G(V_1,E_1)$, \\
$G(V_2,E_2)$ is called a homomorphism of oriented rooted trees \\ $G(V_1,E_1)$, $G(V_2,E_2)$, if $f(\sigma_1)=\sigma_2$, where
 $\sigma_1,\sigma_2$ are the roots of trees $G(V_1,E_1), G(V_2,E_2)$, respectively.
\end{definition}
The homomorphism of oriented rooted trees $f$ is called an {\em isomorphism} of oriented rooted trees, if $f$ is an isomorphism of oriented graphs $G(V_1,E_1), G(V_2,E_2)$. If $G(V_1,E_1)= G(V_2,E_2)$, then homomorphism of oriented rooted trees $f$ is called an {\em endomorphism} of oriented rooted trees. If endomorphism $f$ is an isomorphism, then endomorphism $f$ is called an {\em automorphism} of oriented rooted tree.

\begin{proposition}\label{a5.8.48}
If mappings $f,g$ are automorphisms of oriented rooted tree $G(V,E)$, then $fg$ is an automorphism of tree \\ $G(V,E)$.
\end{proposition}

$\Box$ Suppose that $f,g$ are automorphisms of oriented rooted tree $G(V,E)$ and $s$ is the root of $G(V,E)$. Then we have $sfg=sg=s$. 
If $(u,v)$ is an arc of tree $G(V,E)$, then $(uf,vf)$ also is  an arc of tree $G(V,E)$ (by definitions \ref{d5.8.32} and \ref{d5.8.34} ) because $f$ is an endomorphism. As $(uf,vf)$ is an arc of tree $G(V,E)$, then $(ufg,vfg)$ also is an arc of tree $G(V,E)$ because $g$ also is an endomorphism. Therefore $fg$ is endomorphism of tree $G(V,E)$.

As $f,g$ are automorphisms of tree $G(V,E)$, then $f^{-1}, g^{-1}$ are endomorphisms of the same tree (by definition \ref{d5.8.34}). 
As $f^{-1}, g^{-1}$ are endomorphisms, then $g^{-1}f^{-1}$ is an endomorphism.
Therefore $(fg)^{-1}$ is an endomorphism, because, as $f,g$ are bijections, we have $(fg)^{-1}=g^{-1}f^{-1}$.

Finally, we have proven that $fg$ is an automorphism of oriented rooted tree $G(V,E)$.
\rule{2mm}{2mm}

By ${\rm Aut}(G(V,E))$ we denote set of all automorphisms of oriented rooted tree $G(V,E)$.
\begin{corollary}
$\langle {\rm Aut}(G(V,E)), \otimes \rangle$ is a group. \\
Here
\[
\otimes: \;{\rm Aut}(G(V,E))\times{\rm Aut}(G(V,E)){\to}{\rm Aut}G(V,E)
\]
is the composition of mappings.
\end{corollary}

\section{Machine semigroups}

We would use notation $\mathfrak{M}=\langle Q, A, \circ,*\rangle$ to denote Mealy machine $\mathfrak{M}=\langle Q, A, A, \circ,*\rangle$. In this section the input and output alphabets of all Mealy machines are one and the same, additionally, we drop the requirement that function $Q\times A\stackrel{*}{\to} A$ must be surjective.

For each $q\in Q$ function $Q\times A\stackrel{*}{\to} A$ defines a restricted sequential function
\[
\bar q: A^* \to A^*:u \mapsto q*u.
\]
Fallowing definition is inductive. Let $x\in Q^+$, then with $\overline{xq}$ we denote function 
\[
\overline{xq}: A^* \to A^*:u\mapsto (u\bar x)\bar q\apz \bar q(\bar x(u)).
\]

\begin{definition}
Let $\mathfrak{M}=\langle Q, A, \circ,*\rangle$ be a Mealy machine, then semigroup generated by set $\{\bar q\,|\, q\in Q\}$ is called a machine semigroup $\mathfrak{M}$ (automaton semigroup, semigroup of $\mathfrak{M}$).
\end{definition}
We use notification $\langle \mathfrak{M}\rangle_+$ to denote the semigroup of machine $\mathfrak{M}$.
\begin{corollary}
 $\langle \mathfrak{M} \rangle_+=\{\bar x\,|\, x\in Q^+\}$
\end{corollary} 

Why we are interested in machine semigroups? \\
The definition of machine semigroup implies that it is a convenient
way for represent particular semigroup. In other words if we are given a Mealy machine 
$\mathfrak{M}=\langle Q, A, \circ,*\rangle$, then its machine semigroup is given as well. 

Another question arises: \\
Is there a particular benefit from this representation? How informative is particulat representation?

It turns out that even a very small size Mealy machine gives infinite semigroups.

\begin{examples}  \end{examples}

\begin{figure}[h]
\input{aut45.pic}
\caption{  $|\langle V_{19}\rangle_+|=\infty , \quad  |\langle V_{20} \rangle_+ |=4 $}
\label{z62}
\end{figure}
 
First we prove that semigroup $\langle V_{19}\rangle_+$0 of machine $V_{19}$ is infinite.
We choose infinite word $(01)^\omega\in\overline{0,1}^\omega$, it fallows that
\begin{itemize}
\item $\bar q((01)^\omega)=1^2(10)^\omega$ and $\overline{qp}((01)^\omega)=\bar p(1^2(10)^\omega)=0^2(01)^\omega$
\end{itemize}
Further proof is inductive given  
\[
\overline{(qp)^{n-1}q}((01)^\omega)= 1^{2n}(10)^\omega\quad {\rm and} \quad \overline{(qp)^n}((01)^\omega)=0^{2n}(01)^\omega. 
\]
Hence
\begin{eqnarray*}
\overline{(qp)^n}q((01)^\omega)&=&\bar q(\overline{(qp)^n}((01)^\omega))=\bar q(0^{2n}(01)^\omega) \\
&=&1^{2n}1^2(10)^\omega=1^{2(n+1)}(10)^\omega \\
\rm{ and} \hspace*{35mm}\\
\overline{(qp)^{n+1}}((01)^\omega)&=&\bar p(\overline{q(qp)^n}((01)^\omega)) \\
&=&\bar p (1^{2(n+1)}(10)^\omega)=0^{2(n+1)}(01)^\omega .
\end{eqnarray*}

We have inductively proven that elements of semigroup $\langle V_{19}\rangle_+$
\[
\bar q, \overline{qp}, \overline{qpq}, \ldots, \overline{q(pq)^n}, \overline{(qp)^{n+1}}, \ldots
\]
are unique. Thereby semigroup $\langle V_{19} \rangle_+$ is infinite.

Now lets consider semigroup of machine $V_{20}$. Suppose that $a\in\overline{0,1}$ and $x\in\overline{0,1}^\omega$, then
\begin{itemize}
\item $\bar p(x)=\bar x, \qquad \quad \overline{p^2}(x)=\bar p(\bar p(x))=\bar p(\bar x)=x=\mathbb{I}(x)$
\item $\bar q(ax)=a\bar x, \qquad \overline{q^2}(ax)=\bar q(\bar q(ax))=\bar q(a \bar x)=x=\mathbb{I}(x)$
\item $\overline{qp}(ax)=\bar p(\bar q(ax))=\bar p(\bar ax)=a\bar x, \quad \overline{pq}(ax)=\bar q(\bar p(ax))=q(\bar a\bar x)=\bar ax$ \\
Therefore we have $\overline{pq}=\overline{qp}$, t.i., $\bar p \bar q=\bar q \bar p$.
\end{itemize}
Hence
\[ \begin{array}{c|cccc}
& \mathbb{I} & \bar p & \bar q & \overline{pq}\\
\hline \hline
\mathbb{I} &\mathbb{I}& \bar p & \bar q & \overline{pq}\\
\bar p & \bar p & \mathbb{I} & \overline{pq} & \bar q \\
\bar q & \bar q & \overline{pq} & \mathbb{I} & \bar p\\
\overline{pq} & \overline{pq} & \bar q & \bar p & \mathbb{I} 
\end{array}\]
It fallows that $\langle V_{20} \rangle_+$ is the Klain 4--group $\mathbb{Z}_2\times \mathbb{Z}_2$.

\begin{example}
\textnormal{ Lets choose an alphabet $A=\{a_1,a_2,\ldots, a_p\}$. Lets define a $p$--regular rooted tree $\mathfrak{\bar K}(A)$ as:
\begin{itemize}
\item the words $A^*$ are vertices of the tree;
\item the empty word $\lambda$ is root of the tree;
\item the set $\{(u,ua)\,|\, u\in A^*\wedge a\in A\}$ is the set of arcs.
\end{itemize}}
\end{example}

\begin{figure}[h]
\centering
\scalebox{.8}{\input{aut47.pic}}
\caption{ Labeled rooted trees $\mathfrak{K}$ and $\mathfrak{K}(\overline{0,1})$.}
\label{z70}
\end{figure}

From construction we have that ${\rm deg^-}(\lambda)=0$, $\forall u\in A^+ \; {\rm deg^-}(u)=1$ and $\forall u\in A^* \; {\rm deg^+}(u)=p$. 

\begin{lemma}
The underlying graph of graph $\mathfrak{\bar K}(A)$ is a tree.
\end{lemma}

$\Box$ For the underlying graph of $\mathfrak{\bar K}(A)$ the set of edges is 
$\{[u,ua]\,|\, u\in A^*\wedge a\in A\}$. Lets exclude an arbitrary edge $[u,ua]$. Let prove that for vertex $w$, which reachable from vertex $ua$, holds that $ua\in{\rm Pref} (w)$.
If length of walk is zero, then proposition holds. Suppose that proposition holds for all walks of length $n$ starting at vertex $ua$ .

Lets choose an arbitrary walk 
\[
(ua, l_1, v_1, l_2, v_2,\ldots, l_n,v_n, l_{n+1}, v_{n+1})
\]
of length $n+1$, then $(ua, l_1, v_1, l_2, v_2,\ldots, l_n,v_n)$ is walk of length $n$, and accordingly to the induction hypothesis $ua\in {\rm Pref}(v_n)$.

(i) If $v_n =ua$, then $l_{n+1}=[ua, uab]$ for some $b\in A$ because edge $[u,ua]$ is excluded. Hence
\[
v_{n+1}=uab \quad {\rm and} \quad ua\in {\rm Pref}(uab)={\rm Pref}(v_{n+1}).
\]

(ii) If $v_n \ne ua$, then accordingly with induction hypothesis $v_n=uavb$, where $v\in A^*$ and $b\in A$. Thus $p+1$ edge is incident to vertex $v_n$, i.e.,
\[
[uavb,uav],[uavb,uavba_1],[uavb,uavba_2],\ldots,[uavb,uava_p].
\]
Therefore we have that
\[
v_{n+1}\in\{uav,uavba_1,uavba_2,\ldots,uava_p\}.
\]
Hence $ua\in{\rm Pref}(v_{n+1})$. This concludes inductive part of the proof.

We have proven that vertex $u$ is not reachable. This means that each edge $[u,ua]$ is a bridge, therefore (proposition \ref{a5.8.16}) the underlying graph is a tree.
\rule{2mm}{2mm}
\medskip

As the underlying graph of $\mathfrak{\bar K}(A)$ is a tree, ${\rm deg^-}(\lambda)=0$, $\forall u\in A^+ \; {\rm deg^-}(u)=1$ and  
$\forall u\in A^* \; {\rm deg^+}(u)=p$, thus we have proven, that $\mathfrak{\bar K}(A)$ is a $p$--regular rooted tree. In case when $|A|=2$, the tree $\mathfrak{\bar K}(A)$ is
called a {\em rooted binary tree}.

\begin{proposition}
All rooted $p$--regular trees are isomorphic.
\end{proposition}

$\Box$ Suppose that $G(V,E)$ is a $p$--regular rooted tree with root $\sigma$. Lets prove that $G(V,E)$ is isomorphic to tree $\mathfrak{\bar K}(\overline{0,p-1})$. 

Suppose that $L(n)$ is set of $n$--th level vertices of graph $G(V,E)$. Lets inductively define bijection $f:V\to \overline{0,p-1}^*$. $f( \sigma)\apz \lambda$.

Suppose that $f$ already is defined for all elements of set $\bigcup_{k=0}^n L(k)$, $v\in L(n)$ and $f(v)=u\in\overline{0,p-1}^n$, then ${\rm deg}^+(v)=p$. Suppose that
\[
\centerline{$(v,v_0),(v,v_1),\ldots,(v,v_{p-1})$}
\]
are arcs of tree $G(V,E)$, then
\[
f(v_0)\apz u0, f(v_1)\apz u1, \ldots, f(v_{p-1})\apz u(p-1).
\]
From the definition fallows that $f$ and $f^{-1}$ are bijections. We need to make sure that both mappings are homomorphisms.
Suppose that $(v,w)\in E$, then there exist $n$ such as $v\in L(n)$ and $w\in L(n+1)$.  According to the definition of $f$
\[
\exists u\in\overline{0,p-1}^n,\; f(v)=u \;\wedge\; \exists a\in\overline{0,p-1},\; f(w)=ua.
\]
Therefore
\[
(f(v),f(w))=(u,ua)\in\mathcal{U}\apz\{(u,ua)\,|\, u\in\overline{0,p-1}^* \wedge a\in\overline{0,p-1}\}.
\]
We note that $\mathcal{U}$ is the set of arcs of tree $\mathfrak{\bar K}(\overline{0,p-1})$.

Suppose that $(u,ua)\in\mathcal{U}$, then exists $n$ and $v\in L(n)$ such as $f(v)=u$. According to the definition of $f$ there exist $w\in L(n+1)$ such as $(v,w)\in E$ and $f(w)=ua$.
Hence 
\[
(f^{-1}(u),f^{-1}(ua))=(v,w)\in E.
\]
\rule{2mm}{2mm}

As result with the precision up to isomorphism there exist only one $p$--regular rooted tree. For this reason, unless required otherwise, we choose tree $\mathfrak{\bar K}(\overline{0,p-1})$ as $p$--regular rooted tree.

Let 
\begin{eqnarray}\label{f5.17}
\sigma \stackrel{a_1}{\to} v_1\stackrel{a_2}{\to} \cdots \stackrel{a_n}{\to} v_n
\end{eqnarray}
be a walk in labeled tree $\mathfrak{K}(\overline{0,1})$ (look \ref{z70} fig.), then 
\[
v_1=a_1,v_2=a_1a_2, \ldots, v_n=a_1a_2\ldots a_n.
\] 
Hence if we are interested only in those walks starting at root $\sigma$, we can restrict our attention to $\mathfrak{\bar K}(\overline{0,1})$. Vertex $v_n$ uniquely defines previously mentioned walk (\ref{f5.17}).

We will restrict our attention only on the endomorphisms of tree $\mathfrak{\bar K}(A)$.
With ${\rm End}(A^*)$ we denote the set of all endomorphisms of tree $\mathfrak{\bar K}(A)$.
With ${\rm Aut}(A^*)$ we denote the set of all automorphisms of tree $\mathfrak{\bar K}(A)$.

\begin{lemma}\label{l5.8.36} 
If $f\in{\rm End}(A^*)$, then $\forall u\in A^*\; |f(u)|=|u|$.
\end{lemma}

$\Box$ We have $f(\lambda)=\lambda$ because $\lambda$ is a root. The set or arcs of tree $\mathfrak{\bar K}(A)$ is 
\[
\mathcal{U}=\{(u,ua)\,|\, u\in A^* \wedge a\in A\}.  
\]
As $f$ is an endomorphism, then $\forall a\in A \; f(\lambda,f(a))\in \mathcal{U}$. Hence 
\[
\exists b\in A \; (f(\lambda),f(a))=(\lambda,f(a))=(\lambda,\lambda b)=(\lambda,b).
\]
Therefore $f(a)=b$ and $|f(a)|=|b|=1=|a|$.

Further proof is inductive, given $u\in A^n$, then $|f(u)|=|u|=n$. Suppose that $w\in A^{n+1}$, then $\exists u\in A^n\exists a\in A \; w=ua$.
Hence $(u,w)=(u,ua)\in  \mathcal{U}$.
As $f$ is an  endomorphism, then $(f(u),f(w))=(f(u),f(ua))\in \mathcal{U}$. Hence $|f(ua)|=|f(u)|+1$. Therefore
\[
|f(w)|=|f(ua)|=|f(u)|+1=|u|+1=|ua|=|w|. \quad \rule{2mm}{2mm}
\]

\begin{definition}
$\langle {\rm End}(A^*), \cdot \rangle$ is a monoid. Here 
\[
{\rm End}(A^*)\times{\rm End}(A^*)\stackrel{\cdot}{\to}{\rm End}(A^*)
\]
is a composition of mappings.
\end{definition}

$\Box$ The identity mapping  $\mathbb{I}:A^*\to A^*:u\mapsto u$ is an  endomorphism of tree $\mathfrak{\bar K}(A)$ and serves as neutral element.
We need to prove, that the composition of endomorphisms is an endomorphism.

Suppose that $f,g\in{\rm End}(A^*)$ and $(u,v)$ is an arc of tree $\mathfrak{\bar K}(A)$, then $(f(u),f(v))$ is an arc of tree $\mathfrak{\bar K}(A)$, therefore
$(g(f(u)),g(f(v)))$ is an arc of tree $\mathfrak{\bar K}(A)$.  Finally, we have that $g(f(\lambda))=g(\lambda)=\lambda$. Hence $g\cdot f$ is an endomorphism.
\rule{2mm}{2mm}

\begin{proposition}
A mapping $f:A^*\to A^*$ is a sequential function, if and only if $f\in{\rm End}(A^*)$.
\end{proposition}

$\Box \Rightarrow $ Suppose that  $f:A^*\to A^*$ is sequential function, then $f(\lambda)=\lambda$ because accordingly to the definition of sequential function $|f(\lambda)|=|\lambda|=0$.
Suppose that  $(u,ua)$ is arc of tree $\mathfrak{\bar K}(A)$, then $f(u)\in {\rm Pref}(f(ua))$.  Therefore $|u|=|f(u)|$ and $|ua|=|f(ua)|$. Hence there exist some $b\in A$ such as $f(ua)=f(u)b$. It means that $(f(u),f(ua))=(f(u),f(u)b)$ is arc of tree $\mathfrak{\bar K}(A)$.

$\Leftarrow $ Suppose that $f\in {\rm End}(A^*)$, then (Lemma \ref{l5.8.36}) $\forall u\in A^*\; |f(u)|=|u|$.
Suppose that $u\in A^*$ and $a\in A$, then $(u,ua)$ is arc of tree $\mathfrak{\bar K}(A)$, therefore $(f(u),f(ua))$ is an arc of tree $\mathfrak{\bar K}(A)$. Hence there exist some $b\in A$ such as $f(ua)=f(u)b$. Therefore $f(u)\in {\rm Pref}(f(ua))$. 

Further proof is inductive, given
 \[
u\in \textnormal{Pref}(v)\Rightarrow f(u)\in\textnormal{Pref}(f(v))
\]
only if $|v|-|u|\le n$. 

Suppose that $u\in{\rm Pref}(v)$ and $|v|-|u|=n+1$, then exists $w\in A^*$ and $a\in A$ such as $v=uwa$ and $|w|=n$. As $(uw,uwa)$ is an arc of tree $\mathfrak{\bar K}(A)$, then
$f(uw)\in {\rm Pref}(f(uwa))$. Accordingly with the induction hypothesis $f(u)\in {\rm Pref}(f(uw))$, therefore $f(u)\in {\rm Pref}(f(uwa))={\rm Pref}(f(v))$.
\rule{2mm}{2mm}

\begin{corollary}
If mapping $f:A^*\to A^*$ is a restricted sequential function, then $f\in{\rm End}(A^*)$.
\end{corollary}

Suppose that $q\in Q$, then mapping $\bar q: A^*\to A^*$ sometimes is called $q$ {\em action} on tree $\mathfrak{\bar K}(A)$.

\begin{lemma}
The mapping $\phi: Q^+\to {\rm End}(A^*):x\mapsto \bar x$ is a homomorphism of semigroups.
\end{lemma}

$\Box$ Suppose that $u\in A^*$ and $x,y\in Q^+$, then $u\phi(xy)=u\overline{xy}=(u\bar x)\bar y$ and $(u\phi(x))\phi(y)=(u\bar x)\bar y$. Hence $\phi(xy)=\phi(x)\phi(y)$
\rule{2mm}{2mm}

Suppose that $u\in A^*$ and $x\in Q^+$, then we would use notation $ux\apz u\phi(x)$.
For denoting semigroup ${\rm Im}(\phi)$ we would use notation $\Sigma(\mathfrak{M})$, showing to which machine correspond given semigroup .

\begin{definition}
A semigroup $P$ is called a machine semigroup (automaton semigroup) if there exist a Mealy machine $\mathfrak{M}$ such that $P\cong\Sigma(\mathfrak{M})$. 
\end{definition}

As ${\rm K\overline{er}}(\phi)=\{(x,y)\,|\,\phi(x)=\phi(y)\}$ is a congruence, then $Q^+/\rm{K\overline{er}}$ is semigroup and
$Q^+/\rm{K\overline{er}}(\phi)\cong {\rm Im}(\phi)=\Sigma(\mathfrak{M})$. Accordingly to the $\rm{K\overline{er}}(\phi)$ definition $x$ and $y$ belong to the same coset, if
$\phi(x)=\phi(y)$. Hence
\[
\forall u\in A^* \; ux=u\phi(x)=u\phi(y)=uy.
\]
Therefore $\forall \alpha\in A^\omega\; \alpha x=\alpha y$, where $\alpha x\apz\lim_{n\to \infty} \alpha_n$, $\alpha_n$ is a prefix of $\omega$--word $\alpha$. 

\section{Machine groups}

\begin{definition}
Machine $\mathfrak{M}=\langle Q, A, \circ,*\rangle$ is called to be  invertible if for $\forall q\in Q$ mapping $\bar q|A$ is a bijection. 
\end{definition}

\begin{proposition}
If $\mathfrak{M}=\langle Q, A, \circ,*\rangle$ is an invertible machine, then $\forall q\in Q\; \bar q: A^*\to A^*$ is a bijection.
\end{proposition}

$\Box$ (i) Accordingly to the definition $\bar q|A$ is a bijection. We need to prove that $\forall n\in\mathbb{N}\; \bar q|A^n$ is a bijection.
This in turn would prove that $\bar q:A^*\to A^*$ is a bijection. 

(ii) The proof is done by induction on $n$. Let $w,w'\in A^{n+1}$, then exists words $v, v'\in A^n$ and $a,a'\in A$ such that
\begin{eqnarray*}
w\bar q&=&q*w=q*va=q*v\# q\circ v*a,\\
w'\bar q&=&q*w'=q*v'a'=q*v'\# q\circ v*a'.
\end{eqnarray*} 
Suppose that $\forall q\in Q\; \bar q|A^n$ is a bijection and $v\ne v'$, then $q*v\ne q*v'$, and therefore $w\bar q\ne w'\bar q$.

Suppose that $v=v'$ and $a\ne a'$. 
As for $\forall q\in Q\ $ mapping $\bar q|A^n$ is a bijection, then  $\overline{q\circ v}|A$ is a bijection. Hence $q\circ v*a\ne q\circ v *a'$, and therefore $w\bar q\ne w'\bar q$.
Hence $\bar q|A^{n+1}$ is an injection. 

(iii) Suppose that $w\in A^{n+1}$, then there exist words $v\in A^n$ and $b\in A$ such that $w=vb$. Suppose that $\forall q\in Q\; \bar q|A^n$ is a  bijection, then there exist $u\in A^n$ such that $v= u\bar q=q*u$.

If $\forall q\in Q\; \bar q|A^n$ is a bijection, then $\overline{q\circ u}|A$ is a bijection, therefore exists $a\in A$ such that $q\circ u*a=b$. Hence
\[
ua\bar q=q*ua=q*u\# q\circ u*a=vb=w.
\]
Thus $\bar q|A^{n+1}$ is a surjection.
\rule{2mm}{2mm}

\begin{corollary} \label{s5.8.45} 
If $\mathfrak{M}=\langle Q, A, \circ,*\rangle$ is an invertible machine, then \\
\centerline {$\forall q\in Q\forall v\in A^*\exists ! u\in A^* \; q*u=v$.}
\end{corollary}

$\Box$
$\bar q:A^* \to A^*$ is bijection.
\rule{2mm}{2mm}
\medskip

Lets define for each invertible machine $\mathfrak{M}=\langle Q, A, \circ,*\rangle$  \\
$Q^{-1}\apz\{q^{-1}\,|\, q\in Q\}$;\\
$\bar q^{-1}: A^* \to A^*:v\mapsto vq^{-1}$, where $vq^{-1}\apz u$, if $q*u=v$.

\begin{corollary}\label{s5.8.45}
$\forall q\in Q, \forall u\in A^*, uqq^{-1}=uq^{-1}q=u$.
\end{corollary}

\begin{lemma}\label{l5.8.46} 
If $q\in Q$, then mapping $\bar q^1\in {\rm End}(A^*)$ for each invertible machine $\mathfrak{M}=\langle Q, A, \circ,*\rangle$.
\end{lemma}

$\Box$ The arcs of tree $\mathfrak{K}(A)$ are $(u,ua)$, where $u\in A^*$ and $a\in A$. Mapping $\bar q^1\in {\rm End}(A^*)$ if $(uq^{-1},uaq^{-1})$ is an arc of tree $\mathfrak{K}(A)$.
As we are interested in endomorphisms of rooted oriented tree $\mathfrak{K}(A)$, then we need to prove that $\lambda q^{-1}=\lambda$. Last identity arises from the definition of mapping $\bar q^1$, i.e., 
$\lambda q^{-1}=u'$, if $q*u'=\lambda$. The only word $u'\in A^*$, with property $q*u'=\lambda$, is $u'=\lambda$. Therefore $\lambda q^{-1}=\lambda$.

Suppose that $uq^{-1}=v$ and $uaq^{-1}=w$. According with definition of $\bar q^1$ fallows that,  $q*v=u$ and $q*w=ua$. Suppose that $w=v'b$,  $v'\in A^*$ and $b\in A$, then
\[
ua=q*w=q*v'b=q*v'\#q\circ v'*b.
\]
Hence $q*v'=u=q*v$. As $\bar q$ is a bijection, then $v=v'$. This means that $uaq^{-1}=w=v'b=vb$. Hence $(uq^{-1},uaq^{-1})=(v,w)=(v,vb)$, where $(v,vb)$ is an arc in tree $\mathfrak{K}(A)$.
\rule{2mm}{2mm}
\medskip

For each $q^{-1}\in Q^{-1}$ mapping $\bar q^1\in{\rm End}(A^*)$. We further would use notation $\overline{q^{-1}}\apz \bar q^1$.

Further definition is inductive. Suppose that $x\in (Q\cup Q^{-1})^+$, then with $\overline{xq}$ we would denote function 
\[
\overline{xq}: A^* \to A^*:u\mapsto (u\bar x)\bar q\apz \bar q(\bar x(u)),
\]
if $q\in Q$. Contrary if $q^{-1}\in Q^{-1}$, then with $\overline{xq^{-1}}$ we denote function
\[
\overline{xq^{-1}}: A^* \to A^*:u\mapsto (u\bar x)\bar q^1\apz \bar q^1(\bar x(u)).
\]

\begin{lemma}\label{l5.8.47}
If $\mathfrak{M}=\langle Q, A, \circ,*\rangle$ is an invertible machine, then mapping
\[
\phi: (Q\cup Q^{-1})^+\to {\rm End}(A^*): x\mapsto \bar x
\]
is a homomorphism of semigroups.
\end{lemma}

$\Box$ Suppose that $u\in A^*$ and $x,y\in (Q\cup Q^{-1})^+$, then  $u\phi(xy)=u\overline{xy}=(u\bar x)\bar y$ and $(u\phi(x))\phi(y)=(u\bar x)\bar y$. Therefore $\phi(xy)=\phi(x)\phi(y)$.
\rule{2mm}{2mm}

\begin{definition}
If $\mathfrak{M}=\langle Q, A, \circ,*\rangle$ is invertible machine, then homomorphism
\[
\phi: (Q\cup Q^{-1})^+\to {\rm End}(A^*): x\mapsto \bar x
\]
is called the natural homomorphism.
\end{definition}

\begin{corollary}
If $\mathfrak{M}=\langle Q, A, \circ,*\rangle$ is invertible machine, then \\
$\forall x\in(Q\cup Q^{-1})^+$ mapping $\bar x$ is automorphism of ordered rooted tree $\mathfrak{\bar K}(A)$.
\end{corollary}

$\Box$ Suppose that $q\in Q$.
We have already proved (by lemma \ref{l5.8.46}), that $\bar q^{1}\in{\rm End}(A^*)$. By corollary \ref{s5.8.45} we have that $\bar q\bar q^1=\bar q^1\bar q=\mathbb{I}$, i.e., compositions of those mappings is the identity mapping. Therefore $\bar q^1$ is the inverse mapping of $\bar q$ and an endomorphism. Hence $\bar q,\bar q^1\in{\rm Aut(A^*)}$.

Further proof is inductive. Let $x\in(Q\cup Q^{-1})^+$, $|x|=n$ and $\bar x\in{\rm Aut}(A^*)$. If $q\in Q\cup Q^{-1}$, then  $\bar q\in {\rm Aut(A^*)}$. Hence $\overline{xq}=\bar x \bar q\in{\rm Aut(A^*)}$ (Proposition \ref{a5.8.48}).
\rule{2mm}{2mm}

\begin{proposition}
Let $\mathfrak{M}=\langle Q, A, \circ,*\rangle$ be an invertible machine, then image ${\rm Im}\phi$ of natural homomorphism
\[
\phi: (Q\cup Q^{-1})^+\to {\rm End}(A^*): x\mapsto \bar x
\]
is a subgroup of group ${\rm Aut}(A^*)$.
\end{proposition}

$\Box$
Accordingly to corollary \ref{s5.8.45} we have $\mathbb{I}\in{\rm Im}\phi$. If $f,g\in{\rm Im}\phi$, then exists $x,y\in(Q\cup Q^{-1})^+$ such that $x\phi=f$ and $y\phi=g$.
Hence (by lemma \ref{l5.8.47}) $fg=(x\phi)( y\phi)=xy\phi\in{\rm Im}\phi$.

If $|x|=1$, then $x\in Q\cup Q^{-1}$. Then we have $x=q$ or $x=q^{-1}$ for some $q\in Q$.

(i) Suppose that $x=q$, then $q^{-1}\phi\in{\rm Im}\phi$ and
\[
f(q^{-1}\phi)=(x\phi)(q^{-1}\phi)=(q\phi)(q^{-1}\phi)=(qq^{-1}\phi)=\overline{qq^{-1}}\stackrel{S\ref{s5.8.45}}{=}\mathbb{I}.
\]
Therefore $f^{-1}=q^{-1}\phi\in{\rm Im}\phi$.

(ii) Suppose that $x=q^{-1}$, then $q\phi\in{\rm Im}\phi$ and
\[
fq\phi=(x\phi)(q\phi)=(q^{-1}\phi)(q\phi)=(q^{-1}q)\phi=\overline{q^{-1}q}\stackrel{S\ref{s5.8.45}}{=}\mathbb{I}.
\]
Therefore $f^{-1}=q\phi\in{\rm Im}\phi$.

(iii) Suppose that  $g\in{\rm Im}\phi,\; g=z\phi$, 
$z\in(Q\cup Q^{-1})^k$ and $k\le n$, then exists $z'\in (Q\cup Q^{-1})^k$ such that $g^{-1}=z'\phi\in{\rm Im}\phi$.

Further proof is inductive.
Let $f\in{\rm Im}\phi$ and there exist $x'\in(Q\cup Q^{-1})^{n+1}$ such that $x'\phi=f$.
If so, then here exists $x\in(Q\cup Q^{-1})^n$ and $q\in Q\cup Q^{-1}$ such that $x'=xq$.
Accordingly to the induction hypothesis here exists $y'\in(Q\cup Q^{-1})^n$ and $q'\in Q\cup Q^{-1}$ such that
$(x\phi)^{-1}=y'\phi\in{\rm Im}\phi$ and $(q\phi)^{-1}=q'\phi\in{\rm Im}\phi$. Hence $q'y'\in(Q\cup Q^{-1})^{n+1}$ and $(q'y')\phi\in{\rm Im}\phi$, and we have also
\begin{eqnarray*}
f((q'y')\phi)&=&(x'\phi)((q'y')\phi)=((xq)\phi)((q'y')\phi)=(xqq'y')\phi \\
&=&(x\phi)((qq')\phi)y'\phi=(x\phi)((q\phi)(q'\phi))y'\phi \\
&=&(x\phi)((q\phi)(q\phi)^{-1})y'\phi=(x\phi)\mathbb{I}y'\phi=(x\phi)y'\phi \\
&=&(x\phi)(x\phi)^{-1}=\mathbb{I}.
\end{eqnarray*}
Therefore $f^{-1}=(q'y')\phi\in{\rm Im}\phi$. This concludes the inductive proof.
\rule{2mm}{2mm}

For denoting a group ${\rm Im}\phi$ we would use notation $\Gamma(\mathfrak{M})$ showing the corresponding machine .

\begin{definition}
A group $G$ is called a machine group (automaton group), if there exist an invertible Mealy machine $\mathfrak{M}$ such that $G\cong\Gamma(\mathfrak{M})$. 
\end{definition}

\begin{proposition}
The set
\[
AS_A\apz \{f\in P^A\,|\, f \textnormal{ is bijection}\}
\]
is a group in which group operation is composition of restricted sequential functions.
\end{proposition}

$\Box$
Suppose that $f\in AS_A$, then exists Mealy machine 
\[V=\langle Q, A,A;q_0, \circ, *\rangle,\] such that $\forall u\in A^*\;f(u)=q_0*u$.

Lets define a new Mealy machine $V'=\langle Q, A,A;q_0, \acute \circ, \acute *\rangle$,
where 
\begin{eqnarray} \label{f5.2a1}
 q\acute*a&\apz & b, \textnormal{ ja } q*b=a, \\
\label{f5.2a2} q\acute\circ a&\apz& q\circ b, \textnormal{ ja } q*b=a.
\end{eqnarray}
The new Mealy machine defines restricted sequential function \\
$g(u)\apz q_0\acute* u$.

We need to prove that $g=f^{-1}$, i.e., $g$ is the inverse function of $f$.

(i) Suppose that $Q'\apz\{q\,|\,\exists u\in A^*\;q=q_0\circ u\}$. Lets prove 
\begin{eqnarray}\label{f5.1}
\forall q\in Q'\forall a\in A\; q*(q\acute *a)=a.
\end{eqnarray}

Suppose that $q=q_0\circ u$, then 
\begin{eqnarray*}
f(ua)&=&q_0*ua=q_0*u\#q_0\circ u*a=q_0*u\#q*a.
\end{eqnarray*}
Accordingly to the definition of $q\acute * a$ if $q*b=a$, then $q\acute * a = b$. As $f$ is a bijection, therefore exists unique $b\in A$.
Hence
\begin{eqnarray}\label{f5.2a}
q*(q\acute * a) &=& q*b=a.
\end{eqnarray} 
We have proven, that $fg(a)=a$.

Accordingly to $q\acute \circ a$ we have that
\begin{eqnarray}\label{f5.2b}
q\circ (q\acute * a) &=& q\circ b=q\acute\circ a.
\end{eqnarray} 

Further proof is inductive, given 
\begin{eqnarray}\label{f5.2c}
\forall q'\in Q' \forall v\in A^n \quad q'*(q'\acute * v)=v.
\end{eqnarray}
Suppose that $w\in A^{n+1}$, then there exist $a\in A$ and $v\in A^n$ such that $w=av$. Hence
\begin{eqnarray*}
q*(q\acute * w) &=& q*(q\acute * av) = q*(q\acute * a\# q\acute\circ a\acute* v)\\
&=& q*(q\acute* a)\#q\circ(q\acute*a)*(q\acute\circ a\acute * v)\\
&\underset{(\ref{f5.2a})}{=}&a\#q\circ(q\acute*a)*(q\acute\circ a\acute * v)\\
&\underset{(\ref{f5.2b})}{=}&a\#q\acute\circ a*(q\acute\circ a\acute * v)\\
&\underset{(\ref{f5.2c})}{=}&a\#v=w.
\end{eqnarray*}
This concludes the inductive part of the proof.

In a particular case we have $fg(w)=q_0*(q_0\acute * w)=w$. Therefore $g=f^{-1}$.

(ii) Suppose that $f$ and $g$ are two restricted sequential functions of the set $AS_A$, then there exist two Mealy machines \\
$V=\langle Q, A,A;q_0, \circ, *\rangle$,
 $V'=\langle Q', A,A;q'_0, \acute\circ, \acute*\rangle$ such that
\[
\forall u\in A^* \; (\, f(u)=q_0*u\;\wedge\;g(u)=q'_0\acute* u \,).
\] 
The serial composition of these machines $ V\leadsto V'$ is a machine implementing the composition of functions $f$ and $g$. To be more precise, if 
$\breve V =\langle \breve Q, A,A; \breve q_0, \breve\circ, \breve*\rangle\in  V\leadsto V' $, then 
\[
\forall u\in\overline{0,1}^*\; gf(u)=\breve q_0\breve* u.
\]

\begin{figure}[h] \label{aut35a}
\input{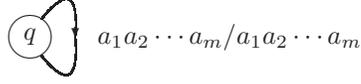}
\caption{Implementation of identity function.}
\end{figure}

(iii) The identity function $\mathbb{I}: A^* \to A^* : u\mapsto u$ is restricted sequential function. Machine implementing identity function given as $V_{15a}$ (shown in \ref{aut35a}fig.).

(iv) We note that the composition of functions is associative.  Finally, we have proven that $AS_A$ is a group.
\rule{2mm}{2mm}

\begin{corollary}
If $\mathfrak{M}=\langle Q, A, \circ,*\rangle$ is an invertible machine, then 
$\Gamma(\mathfrak{M})$ is a subgroup of group $AS_A$.
\end{corollary}

\begin{proposition}
If $|A|=|B|$, then $AS_A\cong AS_B$.
\end{proposition}

$\Box$ If $|A|=|B|$, then exists bijection $\varphi:A \to B$. Suppose that $f\in AS_A$, then exists Mealy machine
\[V=\langle Q, A,A;q_0, \circ, *\rangle,\] ka $\forall u\in A^*\;f(u)=q_0*u$.
Lets define a new Mealy machine $V'=\langle Q, B,B;q_0, \acute \circ, \acute *\rangle$,
where 
\begin{eqnarray*}
q\acute\circ \varphi(a)&\apz& q\circ a, \\
q\acute*\varphi(a)&\apz & \varphi(q*a).
\end{eqnarray*}
Machine $V$ is invertible therefore $\bar q|A:A\to A:a\mapsto q*a$ is a bijection. A composition of bijections is a bijection therefore
\[
\tilde q |B\apz \varphi(\bar q|A)\varphi^{-1}: B\to B: b\mapsto \varphi(q*\varphi^{-1}(b))=q\acute * \varphi\varphi^{-1}(b)=q\acute * b
\]
is a bijection, where $\forall v\in B^*\; \tilde q(v)\apz q\acute * v$. Hence $V'$ is a invertible machine.
This whole construction describes mapping
\[
\psi: AS_A\to AS_B: f\mapsto \psi(f),
\]
where $\psi(f)(v)=q_0\acute* v$ is a restricted sequential function. In present case $f=\bar q_0$ and $\psi(f)=\tilde q_0$.

Suppose that $g\in AS_A$, then there exist a Mealy machine 
\[W=\langle Q', A,A;t_0, \odot, \circledast \rangle,\]
such that $\forall u\in A^*\;g(u)=t_0\circledast u$.
Let define a new Mealy machine\\
$W'\apz \langle Q', A,A;t_0, \acute \odot,\acute \circledast \rangle $,
where 
\begin{eqnarray*}
q\acute\odot \varphi(a)&\apz& q\odot a, \\
q\acute\circledast\varphi(a)&\apz & \varphi(q\circledast a).
\end{eqnarray*}
From previous construction we have that $\tilde t_0=\psi(\bar t_0)=\psi(g)$.
 
The serial composition $ V\leadsto W$ of machines $V, W$ is a machine implementing the composition of functions $f$ and $g$. 
Formally, 
\[
V_W\apz \langle Q'\times Q,A,A;(t_0,q_0),\dot\circ,\dot*\rangle,
\]
where
\begin{eqnarray*}
(t,q)\dot\circ a &\apz& (t \odot q*a,q\circ a);\\
(t,q)\dot*a &\apz& t\circledast q*a,
\end{eqnarray*}
then $\forall u\in A^*\; gf(u)= t_0\circledast q_0*u$.
By previously described construction machine
\[
V'_W\apz \langle Q'\times Q,B,B;(t_0,q_0),\ddot\circ,\ddot*\rangle,
\]
where
\begin{eqnarray*}
(t,q)\ddot\circ \varphi(a) &\apz& (t,q)\dot\circ a, \\
(t,q)\ddot*\varphi(a) &\apz & \varphi((t,q)\ddot*a)=\varphi(t\circledast q* a),
\end{eqnarray*}
implementing mapping $\psi(gf)$, i.e., 
\[
\forall v\in B^* \; \psi(gf)(v)= (t_0,q_0)\ddot * v = \varphi (t_0\circledast q_0 *\varphi^{-1}(v)).
\]
The mapping $\varphi$ is inductively extended on set $A^*$, by a condition
\begin{eqnarray*}
\varphi(ua) &\apz& \varphi(u)\#\varphi(a).
\end{eqnarray*}
\begin{eqnarray*}
(\psi(g)\psi(f))(v) &=& \tilde t_0(\tilde q_0(v))= \tilde t_0(\varphi(q_0*\varphi^{-1}(v)))\\
&=& \varphi (t_0\circledast \varphi^{-1}\varphi(q_0*\varphi^{-1}(v))) \\
&=& \varphi (t_0\circledast (q_0*\varphi^{-1}(v)))=\psi(gf)(v).
\end{eqnarray*}
Therefore mapping $\psi:AS_A \to AS_B$ is a group homomorphism.

Suppose that $f \ne g$, then $\exists u\in A^*\; f(u)\ne g(u)$, i.e.,
\begin{eqnarray*}
q_0*u &=& \bar q_0(u)=f(u)\\
&\ne&g(u) = \bar t_0(u)=t_0\circledast u . 
\end{eqnarray*}
Hence
\begin{eqnarray*}
\psi(f)(\varphi(u)) &=& \tilde q_0(\varphi(u))= q_0\acute* \varphi(u)=\varphi(q_0*u)\\
&\ne&\varphi(t_0\circledast u ) = t_0\acute\circledast \varphi(u)=\tilde t_0(\varphi(u))=\psi(g)(\varphi(u)).
\end{eqnarray*}
Therefore mapping $\psi$ is an injection.

Suppose that $h\in AS_B$, then exists Mealy machine\\
$\check{\mathfrak{M}}=\langle R, B,B;r_0, \check\circ, \check*\rangle$ such that $\forall v\in B^*\; h(v)=r_0\check* v$. We define new machine
$\hat{\mathfrak{M}}\apz\langle R, A,A;r_0, \hat \circ, \hat *\rangle$, where
\begin{eqnarray*}
r\hat\circ \varphi^{-1}(b)&\apz& r\check\circ b, \\
r\hat*\varphi^{-1}(b)&\apz & \varphi^{-1}(r\check*b).
\end{eqnarray*}
For $\hat{\mathfrak{M}}$ we have $\forall u\in A^*\quad \bar r_0(u)\apz r_0\hat* u$. By previously described construction we have that machine $\hat{\mathfrak{M}}$ is invertible, therefore $\bar r_0\in AS_A$.

We apply previously described construction to machine $\hat{\mathfrak{M}}$ and create machine $\breve{\mathfrak{M}}=\langle R, B,B;r_0, \breve\circ, \breve*\rangle$, where
\begin{eqnarray*}
r\breve\circ \varphi(a)&\apz& r\hat\circ a=r\hat\circ\varphi^{-1}\varphi(a)=r\check\circ\varphi(a), \\
r\breve*\varphi(a)&\apz & \varphi(r\hat*a)=\varphi(r\hat*\varphi^{-1}\varphi(a))=\varphi\varphi^{-1}(r\check*\varphi(a))=r\check*\varphi(a).
\end{eqnarray*}
Hence $\psi(\bar r_0)=h$. Therefore $\psi$ is a surjection.

Finally, we have that mapping $\psi$ is an isomorphism of groups $AS_A, AS_B$, therefore $AS_A\cong AS_B$.
\rule{2mm}{2mm}

\begin{example}
\end{example}

\begin{figure}[h]
\input{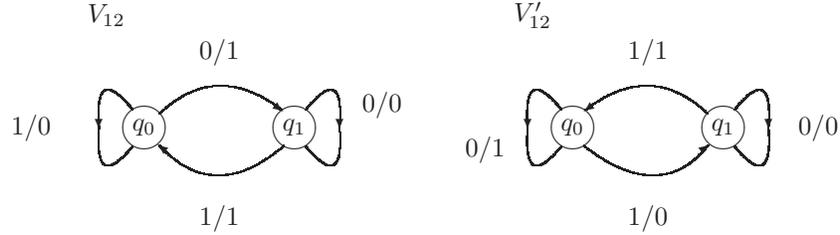}
\caption{ Invertible machine $V_{12}$. }
\label{z2aI}
\end{figure}

Machine $V_{12}$ (\ref{z2aI}fig.) is invertible because $\bar q_0$ is a bijection. Machine $V'_{12}$ (\ref{z2aI}fig.) is constructed according to expressions (\ref{f5.2a1}), (\ref{f5.2a2}):
\begin{eqnarray*}
q_0\acute*0=1, &\textnormal{jo}& q_0*1=0,\\
q_0\acute\circ 0=q_0\circ 1=q_0, &\textnormal{jo}& q_0* 1=0,\\ 
q_0\acute* 1=0, &\textnormal{jo}& q_0* 0=1,\\
q_0\acute\circ 1=q_0\circ 0=q_1, &\textnormal{jo}& q_0*0=1,\\
q_1\acute* 0= 0, &\textnormal{jo}& q_1*0=0,\\
q_1\acute\circ 0= q_1\circ 0=q_1, &\textnormal{jo}& q_1*0=0,\\
q_1\acute* 1= 1, &\textnormal{jo}& q_1*1=1,\\
q_1\acute\circ 1= q_1\circ 1=q_0, &\textnormal{jo}& q_1*1=1.
\end{eqnarray*}
Mapping $\tilde q_0(u)\apz q_0\acute* u$ is inverse of mapping $\bar q_0(u)\apz q_0*u$.
\begin{definition}
Let $\dot V=\langle \dot Q,\dot A, \dot B;\dot\circ,\dot{*} \rangle$, $\ddot V=\langle \ddot Q,\,  \ddot A,\,\ddot B;\ddot\circ, \ddot * \rangle$ be Mealy machines.  A total map 
\[
\mu=(\mu_1,\mu_2,\mu_3):(\dot Q,\dot A, \dot B)\longrightarrow(\ddot Q,\,  \ddot A,\, \ddot B)
\] is called a homomorphism
$\mu :\dot V\longrightarrow\, \ddot V$ if  \\
\[
\mu_1(q\dot\circ a)=\mu_1(q)\ddot\circ\mu_2(a) \;\wedge \;
\mu_3(q\dot*a)=\mu_1(q)\ddot*\mu_2(a)
\]
for  $\forall (q,a)\in Q\times A$.
\end{definition}

Lets consider machines $V_{12}$ u $V'_{12}$ (\ref{z2aI}fig.) and consider mapping
\[
\mu= (\mu_1,\mu_2,\mu_3): V_{12} \to V'_{12},
\]
where
\begin{eqnarray*}
\mu_1 :& q_0 \mapsto q_0,& q_1 \mapsto q_0;\\
\mu_2: &0\mapsto 0,&  1\mapsto 0;\\
\mu_3: &0\mapsto 1,& 1\mapsto 1. 
\end{eqnarray*}
\begin{eqnarray*}
\mu_1(q_0\circ 0) &=& \mu_1(q_1)=q_0=q_0\acute\circ 0=\mu_1(q_0)\acute\circ \mu_2(0),\\
\mu_3(q_0* 0) &=& \mu_3(1)=1=q_0\acute*0=\mu_1(q_0)\acute*\mu_2(0),\\
\mu_1(q_0\circ 1) &=& \mu_1(q_0)=q_0=q_0\acute\circ 0=\mu_1(q_0)\acute\circ \mu_2(1),\\
\mu_3(q_0* 1) &=& \mu_3(0)=1=q_0\acute*0=\mu_1(q_0)\acute*\mu_2(1),\\
\mu_1(q_1\circ 0) &=& \mu_1(q_1)=q_0=q_0\acute\circ 0=\mu_1(q_1)\acute\circ \mu_2(0),\\
\mu_3(q_1* 0) &=& \mu_3(0)=1=q_0\acute*0=\mu_1(q_1)\acute*\mu_2(0),\\
\mu_1(q_1\circ 1) &=& \mu_1(q_0)=q_0=q_0\acute\circ 0=\mu_1(q_1)\acute\circ \mu_2(1),\\
\mu_3(q_1* 1) &=& \mu_3(1)=1=q_0\acute*0=\mu_1(q_1)\acute*\mu_2(1).
\end{eqnarray*}

\begin{proposition}
Let $\dot V=\langle \dot Q, \dot A, \dot B; \dot\circ, \dot* \rangle$ be a Mealy machine, and let there exist state $\dot q\in \dot Q$ and letter of input alphabet $\dot a\in \dot A$ such that
$\dot q\dot\circ\dot a=\dot q$. If so, then for each Mealy machine $\langle Q, A, B;\circ,*\rangle$ there exist a morphism of machines $\mu: V \to \dot V$.
\end{proposition}

$\Box$ Suppose that $\mu=(\mu_1,\mu_2,\mu_3): Q\times A\times B\to \dot Q\times \dot A\times \dot B$ is a mapping defined by conditions 
\begin{eqnarray*}
\forall q\in Q \quad \mu_1(q) &=&\dot q,\\
\forall a\in A \quad \mu_2(a)  &=& \dot a,\\
\forall b\in B \quad \mu_3(b) &=& \dot q\dot*\dot a,
\end{eqnarray*}
then
\begin{eqnarray*}
\mu_1(q\circ a) &=&\dot q = \dot q\dot\circ \dot a=\mu_1(q)\dot\circ\mu_2(a),\\
\mu_3(b) &=& \dot q\dot*\dot a=\mu_1(q)\dot*\mu_2(a).
\end{eqnarray*}
This shows that $\mu: V\to \dot V$ is a homomorphism of machines.
\rule{2mm}{2mm}

In general case this result is not true even for invertible machines. Moreover, as shown in next example, we can choose invertible machines 
\[
V_{13}=\langle Q, A, A;q_0,\circ,*\rangle, V'_{13}=\langle Q', A', A';q'_0,\acute\circ,\acute*\rangle,
\] 
such that mapping $\acute q\acute* u$ is the inverse of mapping $q*u$, but there doesn't exist nether homomorphism $V\to V'$ nor homomorphism $V'\to V$.

\begin{figure}[h]
\input{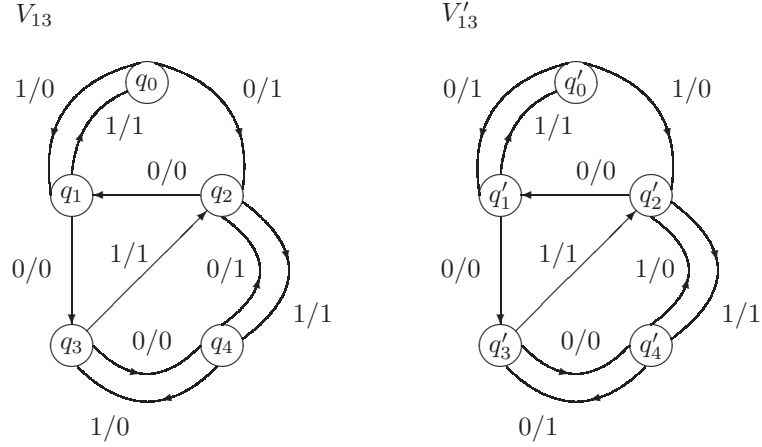}
\caption{ Invertible machine $V_{13}$. }
\label{z2bI}
\end{figure}
 
\begin{example}
\end{example}
Suppose that $\mu=(\mu_1,\mu_2,\mu_3):V_{13} \to V'_{13}$ is a homomorphism.

(i) If $\mu_1(q_0)=q'_0$, then we have cases:
\begin{eqnarray*}
\mu_2: &0\mapsto 0,& 1\mapsto 0;\\
&0\mapsto 0,& 1\mapsto 1;\\
&0\mapsto 1,& 1\mapsto 0;\\
&0\mapsto 1,& 1\mapsto 1.
\end{eqnarray*}

(i{\it a}) $\mu_2: 0\mapsto 0, 1\mapsto 0$.
\begin{eqnarray*}
\mu_1(q_1)&=&\mu_1(q_0\circ 1)=\mu_1(q_0)\acute\circ\mu_2(1)=q'_0\acute\circ 0=q'_1,\\
\mu_3(1)  &=&\mu_3(q_0*0)=\mu_1(q_0)\acute*\mu_2(0)=q'_0\acute* 0=1,\\
\mu_3(1) &=& \mu_3(q_1* 1)=\mu_1(q_1)\acute* \mu_2(1)=q'_1\acute*0=0.
\end{eqnarray*}
Hence $0=\mu_3(1)=1$. A contradiction!

(i{\it b}) $\mu_2: 0\mapsto 0, 1\mapsto 1$.
\begin{eqnarray*}
\mu_1(q_2) &=& \mu_1(q_0\circ 0)= q'_0\acute\circ 0=q'_1,\\
\mu_1(q_1) &=& \mu_1(q_0\circ 1) = q'_0\acute\circ 1=q'_2,\\
\mu_1(q_1) &=& \mu_1(q_2\circ 0)=q'_1\acute\circ 0=q'_3.
\end{eqnarray*}
Hence $q'_3=\mu_1(q_1)=q'_2$. A contradiction!

(i{\it c}) $\mu_2: 0\mapsto 1, 1\mapsto 0$.
\begin{eqnarray*}
\mu_1(q_2) &=& \mu_1(q_0\circ 0)= q'_0\acute\circ 1=q'_2,\\
\mu_1(q_1) &=& \mu_1(q_0\circ 1) = q'_0\acute\circ 0=q'_1,\\
\mu_1(q_1) &=& \mu_1(q_2\circ 0)=q'_2\acute\circ 1=q'_4.
\end{eqnarray*}
Hence $q'_4=\mu_1(q_1)=q'_1$. A contradiction!

(i{\it d}) $\mu_2: 0\mapsto 1, 1\mapsto 1$.
\begin{eqnarray*}
\mu_1(q_2) &=& \mu_1(q_0\circ 0)= q'_0\acute\circ 1=q'_2,\\
\mu_1(q_1) &=& \mu_1(q_0\circ 1) = q'_0\acute\circ 1=q'_2,\\
\mu_1(q_1) &=& \mu_1(q_2\circ 0)=q'_2\acute\circ 1=q'_4.
\end{eqnarray*}
Hence $q'_4=\mu_1(q_1)=q'_2$. A contradiction!

(ii) Suppose that $\mu_1(q_0)=q'_1$.

(ii{\it a}) $\mu_2: 0\mapsto 0, 1\mapsto 0$.
\begin{eqnarray*}
\mu_1(q_2) &=& \mu_1(q_0\circ 0)= q'_1\acute\circ 0=q'_3,\\
\mu_1(q_1) &=& \mu_1(q_0\circ 1) = q'_1\acute\circ 0=q'_3,\\
\mu_1(q_1) &=& \mu_1(q_2\circ 0)=q'_3\acute\circ 0=q'_4.
\end{eqnarray*}
Hence $q'_4=\mu_1(q_1)=q'_3$. A contradiction!

(ii{\it b}) $\mu_2: 0\mapsto 0, 1\mapsto 1$.
\begin{eqnarray*}
\mu_1(q_2) &=& \mu_1(q_0\circ 0)= q'_1\acute\circ 0=q'_3,\\
\mu_1(q_1) &=& \mu_1(q_0\circ 1) = q'_1\acute\circ 1=q'_0,\\
\mu_1(q_1) &=& \mu_1(q_2\circ 0)=q'_3\acute\circ 0=q'_4.
\end{eqnarray*}
Hence $q'_4=\mu_1(q_1)=q'_0$. A contradiction!

(ii{\it c}) $\mu_2: 0\mapsto 1, 1\mapsto 0$.
\begin{eqnarray*}
\mu_1(q_2) &=& \mu_1(q_0\circ 0)= q'_1\acute\circ 1=q'_0,\\
\mu_1(q_1) &=& \mu_1(q_0\circ 1) = q'_1\acute\circ 0=q'_3,\\
\mu_1(q_1) &=& \mu_1(q_2\circ 0)=q'_0\acute\circ 1=q'_2.
\end{eqnarray*}
Hence $q'_2=\mu_1(q_1)=q'_3$. A contradiction!

(ii{\it d}) $\mu_2: 0\mapsto 1, 1\mapsto 1$.
\begin{eqnarray*}
\mu_1(q_2) &=& \mu_1(q_0\circ 0)= q'_1\acute\circ 1=q'_0,\\
\mu_3(0) &=& \mu_3(q_0* 1) = q'_1\acute* 1= 1,\\
\mu_3(0) &=& \mu_3(q_2* 0)=q'_0\acute* 1=0.
\end{eqnarray*}
Hence $0=\mu_3(0)=1$. A contradiction!

(iii) Suppose that $\mu_1(q_0)=q'_2$.

(iii{\it a}) $\mu_2: 0\mapsto 0, 1\mapsto 0$.
\begin{eqnarray*}
\mu_1(q_2) &=& \mu_1(q_0\circ 0)= q'_2\acute\circ 0=q'_1,\\
\mu_1(q_1) &=& \mu_1(q_0\circ 1) = q'_2\acute\circ 0=q'_1,\\
\mu_1(q_1) &=& \mu_1(q_2\circ 0)=q'_1\acute\circ 0=q'_3.
\end{eqnarray*}
Hence $q'_3=\mu_1(q_1)=q'_1$. A contradiction!

(iii{\it b}) $\mu_2: 0\mapsto 0, 1\mapsto 1$.
\begin{eqnarray*}
\mu_1(q_2) &=& \mu_1(q_0\circ 0)= q'_2\acute\circ 0=q'_1,\\
\mu_1(q_1) &=& \mu_1(q_0\circ 1) = q'_2\acute\circ 1=q'_4,\\
\mu_1(q_1) &=& \mu_1(q_2\circ 0)=q'_1\acute\circ 0=q'_3.
\end{eqnarray*}
Hence $q'_3=\mu_1(q_1)=q'_4$. A contradiction!

(iii{\it c}) $\mu_2: 0\mapsto 1, 1\mapsto 0$.
\begin{eqnarray*}
\mu_1(q_2) &=& \mu_1(q_0\circ 0)= q'_2\acute\circ 1=q'_4,\\
\mu_1(q_1) &=& \mu_1(q_0\circ 1) = q'_2\acute\circ 0=q'_1,\\
\mu_1(q_1) &=& \mu_1(q_2\circ 0)=q'_4\acute\circ 0=q'_3.
\end{eqnarray*}
Hence $q'_3=\mu_1(q_1)=q'_1$. A contradiction!
 
(iii{\it d}) $\mu_2: 0\mapsto 1, 1\mapsto 1$.
\begin{eqnarray*}
\mu_1(q_2) &=& \mu_1(q_0\circ 0)= q'_2\acute\circ 1=q'_4,\\
\mu_1(q_1) &=& \mu_1(q_0\circ 1) = q'_2\acute\circ 1=q'_4,\\
\mu_1(q_1) &=& \mu_1(q_2\circ 0)=q'_4\acute\circ 1=q'_3.
\end{eqnarray*}
Hence $q'_3=\mu_1(q_1)=q'_4$. A contradiction!

(iv) Suppose that $\mu_1(q_0)=q'_3$.

(iv{\it a}) $\mu_2: 0\mapsto 0, 1\mapsto 0$.
\begin{eqnarray*}
\mu_1(q_2) &=& \mu_1(q_0\circ 0)= q'_3\acute\circ 0=q'_4,\\
\mu_1(q_1) &=& \mu_1(q_0\circ 1) = q'_3\acute\circ 0=q'_4,\\
\mu_1(q_1) &=& \mu_1(q_2\circ 0)=q'_4\acute\circ 0=q'_3.
\end{eqnarray*}
Hence $q'_3=\mu_1(q_1)=q'_4$. A contradiction!

(iv{\it b}) $\mu_2: 0\mapsto 0, 1\mapsto 1$.
\begin{eqnarray*}
\mu_1(q_2) &=& \mu_1(q_0\circ 0)= q'_3\acute\circ 0=q'_4,\\
\mu_1(q_1) &=& \mu_1(q_0\circ 1) = q'_3\acute\circ 1=q'_2,\\
\mu_1(q_1) &=& \mu_1(q_2\circ 0)=q'_4\acute\circ 0=q'_3.
\end{eqnarray*}
Hence $q'_3=\mu_1(q_1)=q'_2$. A contradiction!

(iv{\it c}) $\mu_2: 0\mapsto 1, 1\mapsto 0$.
\begin{eqnarray*}
\mu_1(q_2) &=& \mu_1(q_0\circ 0)= q'_3\acute\circ 1=q'_2,\\
\mu_3(0) &=& \mu_3(q_0* 1) = q'_3\acute* 0=0,\\
\mu_3(0) &=& \mu_3(q_2* 0)=q'_2\acute* 1=1.
\end{eqnarray*}
Hence $1=\mu_3(0)=0$. A contradiction!
 
(iv{\it d}) $\mu_2: 0\mapsto 1, 1\mapsto 1$.
\begin{eqnarray*}
\mu_1(q_2) &=& \mu_1(q_0\circ 0)= q'_3\acute\circ 1=q'_2,\\
\mu_1(q_1) &=& \mu_1(q_0\circ 1) = q'_3\acute\circ 1=q'_2,\\
\mu_1(q_1) &=& \mu_1(q_2\circ 0)=q'_2\acute\circ 1=q'_4.
\end{eqnarray*}
Hence $q'_4=\mu_1(q_1)=q'_2$. A contradiction!

(v) Suppose that $\mu_1(q_0)=q'_4$.

(v{\it a}) $\mu_2: 0\mapsto 0, 1\mapsto 0$.
\begin{eqnarray*}
\mu_1(q_2) &=& \mu_1(q_0\circ 0)= q'_4\acute\circ 0=q'_3,\\
\mu_1(q_1) &=& \mu_1(q_0\circ 1) = q'_4\acute\circ 0=q'_3,\\
\mu_1(q_1) &=& \mu_1(q_2\circ 0)=q'_3\acute\circ 0=q'_4.
\end{eqnarray*}
Hence $q'_4=\mu_1(q_1)=q'_3$. A contradiction!

(v{\it b}) $\mu_2: 0\mapsto 0, 1\mapsto 1$.
\begin{eqnarray*}
\mu_1(q_2) &=& \mu_1(q_0\circ 0)= q'_4\acute\circ 0=q'_3,\\
\mu_1(q_1) &=& \mu_1(q_0\circ 1) = q'_4\acute\circ 1=q'_2,\\
\mu_1(q_1) &=& \mu_1(q_2\circ 0)=q'_3\acute\circ 0=q'_4.
\end{eqnarray*}
Hence $q'_4=\mu_1(q_1)=q'_2$. A contradiction!

(v{\it c}) $\mu_2: 0\mapsto 1, 1\mapsto 0$.
\begin{eqnarray*}
\mu_1(q_2) &=& \mu_1(q_0\circ 0)= q'_4\acute\circ 1=q'_2,\\
\mu_1(q_1) &=& \mu_1(q_0\circ 1) = q'_4\acute\circ 0=q'_3,\\
\mu_1(q_1) &=& \mu_1(q_2\circ 0)=q'_2\acute\circ 1=q'_4.
\end{eqnarray*}
Hence $q'_4=\mu_1(q_1)=q'_3$. A contradiction!

(v{\it d}) $\mu_2: 0\mapsto 1, 1\mapsto 1$.
\begin{eqnarray*}
\mu_1(q_2) &=& \mu_1(q_0\circ 0)= q'_4\acute\circ 1=q'_2,\\
\mu_1(q_1) &=& \mu_1(q_0\circ 1) = q'_4\acute\circ 1=q'_2,\\
\mu_1(q_1) &=& \mu_1(q_2\circ 0)=q'_2\acute\circ 1=q'_4.
\end{eqnarray*}
Hence $q'_4=\mu_1(q_1)=q'_2$. A contradiction!

In all cases we have arrived at contradictions. This proves that there doesn't exist a homomorphism $\mu=(\mu_1,\mu_2,\mu_3):V_{13} \to V'_{13}$.
\medskip

Suppose that $\mu=(\mu_1,\mu_2,\mu_3):V'_{13} \to V_{13}$ is a homomorphism.

(i) If $\mu_1(q'_0)=q_0$, then we have cases:
\begin{eqnarray*}
\mu_2: &0\mapsto 0,& 1\mapsto 0;\\
&0\mapsto 0,& 1\mapsto 1;\\
&0\mapsto 1,& 1\mapsto 0;\\
&0\mapsto 1,& 1\mapsto 1.
\end{eqnarray*}

(i{\it a}) $\mu_2: 0\mapsto 0, 1\mapsto 0$.
\begin{eqnarray*}
\mu_1(q'_2) &=& \mu_1(q'_0\acute\circ 1)=q_0\circ 0 =q_2,\\
\mu_1(q'_1) &=& \mu_1(q'_2\acute\circ 0)=q_2\circ 0=q_1,\\
\mu_1(q'_1) &=& \mu_1(q'_0\acute\circ 0)=q_0\circ 0=q_2.
\end{eqnarray*}
Hence $q_2=\mu_1(q'_1)=q_1$. A contradiction!

(i{\it b}) $\mu_2: 0\mapsto 0, 1\mapsto 1$.
\begin{eqnarray*}
\mu_1(q'_2) &=& \mu_1(q'_0\acute\circ 1)=q_0\circ 1=q_1,\\
\mu_1(q'_1) &=& \mu_1(q'_2\acute\circ 0)=q_1\circ 0=q_3,\\
\mu_1(q'_1) &=& \mu_1(q'_0\acute\circ 0)=q_0\circ 0=q_2.
\end{eqnarray*}
Hence $q_2=\mu_1(q'_1)=q_3$. A contradiction!

(i{\it c}) $\mu_2: 0\mapsto 1, 1\mapsto 0$.
\begin{eqnarray*}
\mu_1(q'_2) &=& \mu_1(q'_0\acute\circ 1)=q_0\circ 0 =q_2,\\
\mu_1(q'_1) &=& \mu_1(q'_2\acute\circ 0)=q_2\circ 1=q_4,\\
\mu_1(q'_1) &=& \mu_1(q'_0\acute\circ 0)=q_0\circ 1=q_1.
\end{eqnarray*}
Hence $q_1=\mu_1(q'_1)=q_4$. A contradiction!

(i{\it d}) $\mu_2: 0\mapsto 1, 1\mapsto 1$.
\begin{eqnarray*}
\mu_1(q'_2) &=& \mu_1(q'_0\acute\circ 1)=q_0\circ 1 =q_1,\\
\mu_1(q'_1) &=& \mu_1(q'_2\acute\circ 0)=q_1\circ 1=q_0,\\
\mu_1(q'_1) &=& \mu_1(q'_0\acute\circ 0)=q_0\circ 1=q_1.
\end{eqnarray*}
Hence $q_1=\mu_1(q'_1)=q_0$. A contradiction!

(ii) Suppose that $\mu_1(q'_0)=q_1$.

(ii{\it a}) $\mu_2: 0\mapsto 0, 1\mapsto 0$.
\begin{eqnarray*}
\mu_1(q'_2) &=& \mu_1(q'_0\acute\circ 1)=q_1\circ 0 =q_3,\\
\mu_1(q'_1) &=& \mu_1(q'_2\acute\circ 0)=q_3\circ 0=q_4,\\
\mu_1(q'_1) &=& \mu_1(q'_0\acute\circ 0)=q_1\circ 0=q_3.
\end{eqnarray*}
Hence $q_3=\mu_1(q'_1)=q_4$. A contradiction!

(ii{\it b}) $\mu_2: 0\mapsto 0, 1\mapsto 1$.
\begin{eqnarray*}
\mu_1(q'_2) &=& \mu_1(q'_0\acute\circ 1)=q_1\circ 1=q_0,\\
\mu_1(q'_1) &=& \mu_1(q'_2\acute\circ 0)=q_0\circ 0=q_2,\\
\mu_1(q'_1) &=& \mu_1(q'_0\acute\circ 0)=q_1\circ 0=q_3.
\end{eqnarray*}
Hence $q_3=\mu_1(q'_1)=q_2$. A contradiction!

(ii{\it c}) $\mu_2: 0\mapsto 1, 1\mapsto 0$.
\begin{eqnarray*}
\mu_1(q'_2) &=& \mu_1(q'_0\acute\circ 1)=q_1\circ 0 =q_3,\\
\mu_1(q'_1) &=& \mu_1(q'_2\acute\circ 0)=q_3\circ 1=q_4,\\
\mu_1(q'_1) &=& \mu_1(q'_0\acute\circ 0)=q_1\circ 1=q_0.
\end{eqnarray*}
Hence $q_0=\mu_1(q'_1)=q_4$. A contradiction!

(ii{\it d}) $\mu_2: 0\mapsto 1, 1\mapsto 1$.
\begin{eqnarray*}
\mu_1(q'_2) &=& \mu_1(q'_0\acute\circ 1)=q_1\circ 1 =q_0,\\
\mu_1(q'_1) &=& \mu_1(q'_2\acute\circ 0)=q_0\circ 1=q_1,\\
\mu_1(q'_1) &=& \mu_1(q'_0\acute\circ 0)=q_1\circ 1=q_0.
\end{eqnarray*}
Hence $q_0=\mu_1(q'_1)=q_1$. A contradiction!

(iii) Suppose that $\mu_1(q'_0)=q_2$.

(iii{\it a}) $\mu_2: 0\mapsto 0, 1\mapsto 0$.
\begin{eqnarray*}
\mu_1(q'_2) &=& \mu_1(q'_0\acute\circ 1)=q_2\circ 0 =q_1,\\
\mu_1(q'_1) &=& \mu_1(q'_2\acute\circ 0)=q_1\circ 0=q_3,\\
\mu_1(q'_1) &=& \mu_1(q'_0\acute\circ 0)=q_2\circ 0=q_1.
\end{eqnarray*}
Hence $q_1=\mu_1(q'_1)=q_3$. A contradiction!

(iii{\it b}) $\mu_2: 0\mapsto 0, 1\mapsto 1$.
\begin{eqnarray*}
\mu_1(q'_2) &=& \mu_1(q'_0\acute\circ 1)=q_2\circ 1=q_4,\\
\mu_1(q'_1) &=& \mu_1(q'_2\acute\circ 0)=q_4\circ 0=q_2,\\
\mu_1(q'_1) &=& \mu_1(q'_0\acute\circ 0)=q_2\circ 0=q_1.
\end{eqnarray*}
Hence $q_1=\mu_1(q'_1)=q_2$. A contradiction!

(iii{\it c}) $\mu_2: 0\mapsto 1, 1\mapsto 0$.
\begin{eqnarray*}
\mu_1(q'_2) &=& \mu_1(q'_0\acute\circ 1)=q_2\circ 0 =q_1,\\
\mu_1(q'_1) &=& \mu_1(q'_2\acute\circ 0)=q_1\circ 1=q_0,\\
\mu_1(q'_1) &=& \mu_1(q'_0\acute\circ 0)=q_2\circ 1=q_4.
\end{eqnarray*}
Hence $q_4=\mu_1(q'_1)=q_0$. A contradiction!

(iii{\it d}) $\mu_2: 0\mapsto 1, 1\mapsto 1$.
\begin{eqnarray*}
\mu_1(q'_2) &=& \mu_1(q'_0\acute\circ 1)=q_2\circ 1 =q_4,\\
\mu_1(q'_1) &=& \mu_1(q'_2\acute\circ 0)=q_4\circ 1=q_3,\\
\mu_1(q'_1) &=& \mu_1(q'_0\acute\circ 0)=q_2\circ 1=q_4.
\end{eqnarray*}

(iv) Suppose that $\mu_1(q'_0)=q_3$.

(iv{\it a}) $\mu_2: 0\mapsto 0, 1\mapsto 0$.
\begin{eqnarray*}
\mu_1(q'_2) &=& \mu_1(q'_0\acute\circ 1)=q_3\circ 0 =q_4,\\
\mu_1(q'_1) &=& \mu_1(q'_2\acute\circ 0)=q_4\circ 0=q_2,\\
\mu_1(q'_1) &=& \mu_1(q'_0\acute\circ 0)=q_3\circ 0=q_4.
\end{eqnarray*}
Hence $q_4=\mu_1(q'_1)=q_2$. A contradiction!

(iv{\it b}) $\mu_2: 0\mapsto 0, 1\mapsto 1$.
\begin{eqnarray*}
\mu_1(q'_2) &=& \mu_1(q'_0\acute\circ 1)=q_3\circ 1=q_2,\\
\mu_1(q'_1) &=& \mu_1(q'_2\acute\circ 0)=q_2\circ 0=q_1,\\
\mu_1(q'_1) &=& \mu_1(q'_0\acute\circ 0)=q_3\circ 0=q_4.
\end{eqnarray*}
Hence $q_4=\mu_1(q'_1)=q_3$. A contradiction!

(iv{\it c}) $\mu_2: 0\mapsto 1, 1\mapsto 0$.
\begin{eqnarray*}
\mu_1(q'_2) &=& \mu_1(q'_0\acute\circ 1)=q_3\circ 0 =q_4,\\
\mu_1(q'_1) &=& \mu_1(q'_2\acute\circ 0)=q_4\circ 1=q_3,\\
\mu_1(q'_1) &=& \mu_1(q'_0\acute\circ 0)=q_3\circ 1=q_2.
\end{eqnarray*}
Hence $q_2=\mu_1(q'_1)=q_3$. A contradiction!

(iv{\it d}) $\mu_2: 0\mapsto 1, 1\mapsto 1$.
\begin{eqnarray*}
\mu_1(q'_2) &=& \mu_1(q'_0\acute\circ 1)=q_3\circ 1 =q_2,\\
\mu_1(q'_1) &=& \mu_1(q'_2\acute\circ 0)=q_2\circ 1=q_4,\\
\mu_1(q'_1) &=& \mu_1(q'_0\acute\circ 0)=q_3\circ 1=q_2.
\end{eqnarray*}
Hence $q_2=\mu_1(q'_1)=q_4$. A contradiction!

(v) Suppose that $\mu_1(q'_0)=q_4$.

(v{\it a}) $\mu_2: 0\mapsto 0, 1\mapsto 0$.
\begin{eqnarray*}
\mu_1(q'_2) &=& \mu_1(q'_0\acute\circ 1)=q_4\circ 0 =q_2,\\
\mu_1(q'_1) &=& \mu_1(q'_2\acute\circ 0)=q_2\circ 0=q_1,\\
\mu_1(q'_1) &=& \mu_1(q'_0\acute\circ 0)=q_4\circ 0=q_2.
\end{eqnarray*}
Hence $q_2=\mu_1(q'_1)=q_1$. A contradiction!

(v{\it b}) $\mu_2: 0\mapsto 0, 1\mapsto 1$.
\begin{eqnarray*}
\mu_1(q'_2) &=& \mu_1(q'_0\acute\circ 1)=q_4\circ 1=q_3,\\
\mu_1(q'_1) &=& \mu_1(q'_2\acute\circ 0)=q_3\circ 0=q_4,\\
\mu_1(q'_1) &=& \mu_1(q'_0\acute\circ 0)=q_4\circ 0=q_2.
\end{eqnarray*}
Therefore $q_2=\mu_1(q'_1)=q_4$. A contradiction!

(v{\it c}) $\mu_2: 0\mapsto 1, 1\mapsto 0$.
\begin{eqnarray*}
\mu_1(q'_2) &=& \mu_1(q'_0\acute\circ 1)=q_4\circ 0 =q_2,\\
\mu_1(q'_1) &=& \mu_1(q'_2\acute\circ 0)=q_2\circ 1=q_0,\\
\mu_1(q'_1) &=& \mu_1(q'_0\acute\circ 0)=q_4\circ 1=q_3.
\end{eqnarray*}
Therefore $q_3=\mu_1(q'_1)=q_0$. A contradiction!

(v{\it d}) $\mu_2: 0\mapsto 1, 1\mapsto 1$.
\begin{eqnarray*}
\mu_1(q'_2) &=& \mu_1(q'_0\acute\circ 1)=q_4\circ 1 =q_3,\\
\mu_1(q'_1) &=& \mu_1(q'_2\acute\circ 0)=q_3\circ 1=q_2,\\
\mu_1(q'_1) &=& \mu_1(q'_0\acute\circ 0)=q_4\circ 1=q_3.
\end{eqnarray*}
Therefore $q_3=\mu_1(q'_1)=q_2$. A contradiction!

In all cases we have arrived at contradictions. This proves that there doesn't exist a homomorphism $\mu=(\mu_1,\mu_2,\mu_3):V'_{13} \to V_{13}$.

\begin{definition}
We say that machine $'V$ simulates machine $V$ if there exist
mappings
\[
Q\stackrel{h_1}{\longrightarrow}{'Q},\quad 
A\stackrel{h_2}{\longrightarrow}{'\!\!A},\quad
'\!B\stackrel{h_3}{\longrightarrow}{B}
\]
such that  the diagram
\[
\begin{array}{cccccc}
&Q & \times & {A^*} & \stackrel{\ast}{\longrightarrow} & {B^*}\\
\lefteqn{h_1}&
\downarrow && \downarrow \lefteqn{h_2} && \uparrow \lefteqn{h_3}\\
& 'Q & \times & {'\!\!A^*} & \stackrel{\ast}{\longrightarrow} & {'\!B^*}
\end{array}
\]
commutes. That is if $\forall(q,{u})\in {Q}\times {A^*} \; q\ast {u}=h_3(h_1({q})\ast h_2(u))$.
\end{definition}

As it turns out (consider next example), we can choose such invertible machines \[V_{12}=\langle Q, A, A;q_0,\circ,*\rangle, V'_{12}=\langle Q', A', A';q_0,\acute\circ,\acute*\rangle,\]
that $\acute q\acute* u$ is the inverse mapping of $q*u$, but machine $V'_{12}$ doesn't simulate machine $V_{12}$, and $V_{12}$ doesn't simulate $V'_{12}$.

\begin{example}
\end{example}

Suppose that machine $V'_{12}$ simulates machine $V_{12}$ (\ref{z2aI} fig.), i.e., exists mappings $h_1,h_2,h_3$, that
\[
q*u=h_3(h_2(q)\acute*h_2(u))
\]
for $\forall u\in\overline{0,1}^*$ and for all states $q$ of machine $V_{12}$. Therefore
\begin{eqnarray*}
0 &=& q_0*1=h_3(h_1(q_0))\acute*h_2(1)),\\
1 &=& q_0*0=h_3(h_1(q_0)\acute* h_2(0)).
\end{eqnarray*}
This shows that mappings $h_2$ and $h_3$ are bijections. Therefore exists inverse mapping (a bijection) $h$ of mapping $h_3$. Hence
\[
h(q*u)=h(h_3(h_1(q)\acute* h_2(u)))=h_1(q)\acute*h_2(u).
\]
We also have that
\begin{eqnarray*}
h(1)=h(q_0*0) &=& h_1(q_0)\acute* h_2(0),\\
h(0)=h(q_1*0) &=& h_1(q_1)\acute* h_2(0).
\end{eqnarray*}
As $h$ is bijection, then we have that $h_1$ also is bijection. Hence all mappings $h,h_1,h_2$ are bijections.

(i) Suppose that $h_1: q_0\mapsto q_0, q_1\mapsto q_1$, then
\[
h(10)=h(q_0*00)=h_1(q_0)\acute* h_2(00))=q_0\acute* h_2(00).
\] 
As $h$ is bijection, then $h_2(0)=1$. Mapping $h_2$ also is a bijection, therefore $h_2(1)=0$. Hence
\[
h(10)=h(q_1*11)=h_1(q_1)\acute*h_2(11)=q_1\acute* 00=00.
\]
Contradiction because $h$ is a bijection.

(ii) Suppose that $h_1: q_0\mapsto q_1, q_1\mapsto q_0$, then
\[
h(10)=h(q_0*00)=h_1(q_0)\acute* h_2(00))=q_1\acute* h_2(00).
\]
As $h$ is a bijection, then $h_2(0)=1$. Mapping $h_2$ also is a bijection, therefore $h_2(1)=0$. Hence
\[
h(10)=h(q_1*11)=h_2(q_1)\acute*h_2(11)=q_0\acute*00=11.
\]
This is a contradiction because $h$ is a bijection.

We have obtained contradictions in both cases, consequentially machine  $V'_{12}$ is not capable of simulating machine $V_{12}$.

Suppose that machine $V_{12}$ is simulating machine $V'_{12}$ (look at \ref{z2aI} fig.), i.e., there exist mappings $h_1,h_2,h_3$ such that
\[
q\acute*u=h_3(h_2(q)*h_2(u))
\]
for all $u\in\overline{0,1}^*$ and for all states $q$ of machine $V'_{12}$. Therefore
\begin{eqnarray*}
0 &=& q_0\acute*1=h_3(h_1(q_0))*h_2(1)),\\
1 &=& q_0\acute*0=h_3(h_1(q_0)* h_2(0)).
\end{eqnarray*}
This shows that $h_2$ and $h_3$ are bijections. Hence there exist a inverse mapping $h$ of mapping $h_3$, which also is a bijection. Hence
\[
h(q\acute*u)=h(h_3(h_1(q)* h_2(u)))=h_1(q)*h_2(u).
\]
We also have that
\begin{eqnarray*}
h(1)=h(q_0\acute*0) &=& h_1(q_0)* h_2(0),\\
h(0)=h(q_1\acute*0) &=& h_1(q_1)* h_2(0).
\end{eqnarray*}
Therefore $h$ is a bijection, consequentially we have, that $h_1$ is a bijection. Hence all mappings $h,h_1,h_2$ are bijections.

(i) Suppose that $h_1: q_0\mapsto q_0, q_1\mapsto q_1$, then
\[
h(01)=h(q_0\acute*11)=h_1(q_0)* h_2(11)=q_0* h_2(11).
\] 
As $h$ is a bijection, then $h_2(1)=0$. Hence
\[
h(10)=h(q_1\acute*11)=h_1(q_1)*h_2(11)=q_1* 00=00.
\]
This is a contradiction because $h$ is a bijection.

(ii) Suppose that $h_1: q_0\mapsto q_1, q_1\mapsto q_0$, then
\[
h(10)=h(q_1\acute*11)=h_1(q_1)* h_2(11))=q_0* h_2(11).
\]
As $h$ is a bijection, then $h_2(1)=0$. The mapping $h_2$ also is a bijection, therefore $h_2(1)=0$. Hence
\[
h(01)=h(q_0\acute*11)=h_2(q_0)*h_2(11)=q_1*00=00.
\]
This is a contradiction because $h$ is a bijection.

Thus we have obtained contradictions in both cases, consequentially the machine $V_{12}$ is not capable of simulating $V'_{12}$.


\begin{thebibliography}{99}

\bibitem{kudr} {\cyr Kudryavcev V. B.,
 Aleshin S. V., Podkolzin A. S.} (1985) \linebreak 
{\cyi Vvedenie v teoriyu avtomatov.} [\,{\em An Introduction to the 
\linebreak Theory 
of Automata.}\,] 
{\cyr Moskva <Nauka>.} (Russian)

\bibitem{ruoh}
{Ruohonen K. (2008)} {\em Graph Theory.} \\
{\color{blue} \underline {https://archive.org/details/flooved3467}}
{(accessed 26. december, 2017)}

\bibitem{Boll}
{Bollobas B. (1998)} {\em Modern Graph Theory.} {Springer--Verlag.}

\bibitem{Cohn}
{Cohn P.M. (1981)} {\em Universal Algebra.} {Springer--Verlag.}

\bibitem{Lang} 
{Lang S. (2005)} {\em Undergraduate Algebra.} {Springer--Verlag.}

\bibitem{Pin}
{Pin J.E. (2016)} {\em Mathematical Foundations of Automata Theory.}
{\color{blue} \underline {https://www.irif.fr/~jep/PDF/MPRI/MPRI.pdf}}
{(accessed 14. january, 2018)} 
\end{thebibliography}
\end{document}